\documentclass[11pt]{article}
\usepackage{epsfig}
\usepackage{amssymb,amsmath,amsthm,amscd}
\usepackage{latexsym}

\newtheorem{thm}{Theorem}[section]
\newtheorem{prop}[thm]{Proposition}

\newtheorem{lem}[thm]{Lemma}

\theoremstyle{definition}
\newtheorem{defn}[thm]{Definition}

\newcounter{labelflag} 

\newcommand{\Label}[1]{
                       \ifnum\thelabelflag=1
                          \ifmmode
                             \makebox[0in][l]{\qquad\fbox{\rm#1}}
                          \else
                             \marginpar{\vspace{0.7\baselineskip}
                                        \hspace{-1.1\textwidth}
                                        \fbox{\rm#1}}
                          \fi
                       \fi
                       \label{#1}
                      }

\newcommand{\be}{\begin{equation}}
\newcommand{\ee}{\end{equation}}

\newcommand{\h}{L^2(Q)}
\newcommand{\hone}{H^1_0(Q)}

\newcommand{\rone}{\mathbb R}
\newcommand{\rtwo}{\mathbb R^2}

\newcommand{\De}{\Delta}

\newcommand{\rn}{{\mathbf R}^n}

 \newcommand{\R}{\mathbb{R}}

\newcommand{\bbme}{Benjamin-Bona-Mahony equation }
\newcommand{\bbm}{Benjamin-Bona-Mahony  }
\newcommand{\vecf}{\overrightarrow{F}}
\newcommand{\vecg}{\overrightarrow{G}}
\newcommand{\zt}{z(\theta_t \omega)}
\newcommand{\cutf}{\phi^{2}({\frac{x_{3}^{2}}{k^{2}}})}
\newcommand{\vb}{{\tilde{v}}}
\begin{document}

\begin{titlepage}
\title{\Large\bf   Random Attractors for  the  Stochastic
 Benjamin-Bona-Mahony Equation on Unbounded Domains}
\vspace{7mm}

\author{
Bixiang Wang  \thanks {Supported in part by NSF  grant DMS-0703521}
\vspace{1mm}\\
Department of Mathematics, New Mexico Institute of Mining and
Technology \vspace{1mm}\\ Socorro,  NM~87801, USA \vspace{1mm}\\
Email: bwang@nmt.edu}
\date{}
\end{titlepage}

\maketitle

\medskip

\begin{abstract}
 We prove the  existence of
 a compact random attractor for the
stochastic  Benjamin-Bona-Mahony Equation
defined on an unbounded domain.
This  random  attractor is invariant and attracts
 every pulled-back tempered random set
 under the forward flow.
The asymptotic compactness of the  random dynamical system
   is established by a
   tail-estimates method,    which shows that the solutions are 
   uniformly 
   asymptotically  small when space and time variables approach
   infinity.

\end{abstract}

{\bf Key words.}    Stochastic  Benjamin-Bona-Mahony Equation,    random attractor, pullback attractor,
asymptotic compactness.

 {\bf MSC 2000.} Primary 60H15.  Secondary  35B40,  35B41.

%%%%%%%%%%%%%%%%%%%%%%%%%%%%%%%%%%%%%%%%%%%%%%%%%%%%%%%%
%  Running title, etc.
%%%%%%%%%%%%%%%%%%%%%%%%%%%%%%%%%%%%%%%%%%%%%%%%%%%%%%%%
%

%%%%%%%%%%%%%%%%%%%%%%%%%%%%%%%%%%%%%%%%%%%%%%%%%%%%%%%%
%  Body of article
%%%%%%%%%%%%%%%%%%%%%%%%%%%%%%%%%%%%%%%%%%%%%%%%%%%%%%%%
%
 %\renewcommand{\baselinestretch}{1.5}
%

\baselineskip=1.5\baselineskip

\section{Introduction}
\setcounter{equation}{0}
This paper is concerned with  the asymptotic behavior of solutions of the
stochastic 
\bbm (BBM) equation  on  an unbounded  three-dimensional channel.
Let $Q=D \times \rone$ where $D$ is a
bounded open subset of $\rtwo$. Consider the following   BBM  equation   on $Q$:
\be
\label{intro1}
d u -d (\De u) - \nu \De u dt + \nabla \cdot {\overrightarrow{F}}(u)  dt
= g dt + h  dw,   \quad x \in Q, \quad t>0,
\ee
 where $\nu$ is a positive constant, $g = g(x)$ and $h=h(x)$ are given functions
defined on $Q$,  
  $\vecf $   is a  smooth  nonlinear vector function, 
  and
  $w$   is a    two-sided real-valued
 Wiener process  on a probability
 space which will be specified later.

Stochastic differential equations arise from many physical
 systems when random uncertainties are taken into account.
 The long-term behavior of   random systems is captured
by a pullback random attractor. This concept
 was introduced in \cite{cra2, fla1} as
  an  extension to    stochastic  systems
 of  the  attractors theory  of
  deterministic   equations in \cite{bab1, hal1, rob1, sel1, tem1}.
  The existence of   random attractors has been studied
  extensively by many authors, see, e.g.,  \cite{arn1, bat1, car1, cra1, cra2, fla1}
and the references therein.
Notice that
    the partial  differential  equations (PDEs)   studied in these papers
  are all defined in {\it bounded} domains.
  In the case of {\it  unbounded}  domains, the existence of random attractors
  for PDEs
  is not  well  understood yet. 
  In this case, as far as we know,  the existence of random attractors
  is only established  for the Reaction-Diffusion equation on $\R^n$
 in \cite{bat2}
  recently.
 In this paper, we will     investigate 
   the  existence of    a  compact  random attractor  for the stochastic BBM equation
 defined on the {\it unbounded} channel   $Q$.

Notice that 
  Sobolev embeddings are no longer compact 
for the unbounded domain $Q$.
This
introduces a major obstacle  for proving the existence of   attractors
for the \bbme  defined on $Q$.
  For some  deterministic equations,  the difficulty caused by
the unboundedness of domains  may  be overcome
by  the
  energy equation approach or by the  tail-estimates approach.
The energy equation method  was   developed by  Ball in
\cite{bal1, bal2}  and used by many authors  (see, e.g.,   
  \cite{ghi1,  gou1, ju1, moi2, moi1,  ros1, wanx}).
  The
   tail-estimates approach      was  developed  in \cite{wan1}
   for  {\it deterministic}   PDEs
   and used in
    \cite{ant1, ant2, arr1,  mor1, pri1, rod1, sta1, sun1}.
    In this paper, we will develop a tail-estimates approach for 
     weakly dissipative
{\it stochastic}  PDEs like  the  \bbme    and  prove the
    existence of   compact random attractors  on unbounded domains.  The idea is based on the 
    observation 
    that the solutions of the equation are uniformly small when space and time
    variables are sufficiently large.

We mention that the \bbme  
  was proposed  in \cite{ben-bon-mah} as  a
 model for propagation of long waves which
  incorporates nonlinear  dispersive and
dissipative effects.
In the deterministic case, 
the  existence  and  uniqueness  of solutions
  were  studied  in
 \cite{avr, avr-gol, ben-bon-mah,
bon-bry, bon-dou, che, gol-wic, med-mil, med-per},
and the global attractors were investigated in    \cite{ast-bis-bis-fer, cel-kal-pol,
chu-pol-sie, sta1,
wan2,
wan3, wan-yan}.
In order  to deal with the stochastic \bbm equation, we need to
transform the stochastic equation with a random term into
a deterministic one with a random parameter. This transformation
will change the structure of the original equation and hence cause
some extra difficulties for  deriving uniform estimates on solutions,
especially  on the tails of solutions for large space and time variables.
For instance, if we   take the inner product of   the deterministic 
BBM equation  with
$u$ in $L^2(Q)$, then the nonlinear term disappears, and hence 
the uniform estimates   in this case are not hard to get.
However, after transformation, this property is lost and the 
nonlinear term does not disappear when performing energy
estimates. This is the reason why much effort of this paper
is devoted to deriving the uniform estimates on the tails of
solutions (see Section 4 for details).

 This paper is organized  as follows. In the next section, we
 review the  pullback
 random attractors  theory  for random dynamical systems. In Section 3,
 we define  a continuous  random dynamical system for the
 stochastic \bbme  on $Q$.
 Then we derive the uniform estimates of solutions  in Section 4, which include
 the uniform estimates on the tails of solutions. 
Finally,  in Section 5, we   establish the asymptotic compactness of
the random dynamical system and prove
  the existence of   a  pullback random attractor.

In the sequel,  we adopt  the following notations.  We denote by
$\| \cdot \|$ and $(\cdot, \cdot)$ the norm and the inner product
of
  $L^2(Q)$, respectively.    The
norm of  a given  Banach space $X$  is written as    $\|\cdot\|_{X}$.
We also use $\| \cdot\|_{p}$    to denote   the norm  of
$L^{p}(Q)$.  The letters $c$ and $c_i$ ($i=1, 2, \ldots$)
are  generic positive constants  which may change their  values from line to
line or even in the same line.

Throughout this paper, we will frequently use the embedding inequality
\be
\label{agmon}
\| u \|_{\infty} \le \beta_0 \| u \|_{H^2(Q)}, \quad \forall \ u \in H^2(Q),
\ee
and the Poincare inequality
\be
\label{poincare}
\|  \nabla u \| ^2  \ge \lambda \| u \| ^2  , \quad  \forall \ u \in   H^1_0(Q),
\ee
where $\beta_0$ and $\lambda$ are positive constants.

\section{Preliminaries}
\setcounter{equation}{0}

In this section, we recall some basic concepts
related to random attractors for stochastic dynamical
systems. The reader is referred to \cite{arn1, bat1, cra1, fla1} for more details.

Let  $(X, \| \cdot \|_X)$ be a   separable
Hilbert space with Borel $\sigma$-algebra $\mathcal{B}(X)$,
 and
$(\Omega, \mathcal{F}, P)$  be  a probability space.

\begin{defn}
$(\Omega, \mathcal{F}, P,  (\theta_t)_{t\in \R})$
is called a  metric   dynamical  system
if $\theta: \R \times \ \Omega \to \Omega$ is
$(\mathcal{B}(\R) \times \mathcal{F}, \mathcal{F})$-measurable,
$\theta_0$ is the identity on $\Omega$,
$\theta_{s+t} = \theta_t \circ  \theta_s$ for all
$s, t \in \R$ and $\theta_t P = P$ for all $t \in \R$.
 \end{defn}

\begin{defn}
\label{RDS}
A continuous random dynamical system (RDS)
on $ X  $ over  a metric  dynamical system
$(\Omega, \mathcal{F}, P,  (\theta_t)_{t\in \R})$
is  a mapping
  $$
\Phi: \R^+ \times \Omega \times X \to X \quad (t, \omega, x)
\mapsto \Phi(t, \omega, x),
$$
which is $(\mathcal{B}(\R^+) \times \mathcal{F} \times \mathcal{B}(X), \mathcal{B}(X))$-measurable and
satisfies, for $P$-a.e.  $\omega \in \Omega$,

(i) \  $\Phi(0, \omega, \cdot) $ is the identity on $X$;

(ii) \  $\Phi(t+s, \omega, \cdot) = \Phi(t, \theta_s \omega,
\cdot) \circ \Phi(s, \omega, \cdot)$ for all $t, s \in \R^+$;

(iii) \  $\Phi(t, \omega, \cdot): X \to  X$ is continuous for all
$t \in  \R^+$.
\end{defn}

Hereafter, we always assume that $\Phi$  is a continuous RDS on $X$
over $(\Omega, \mathcal{F}, P,  (\theta_t)_{t\in \R})$.

\begin{defn}
 A random  bounded set $\{B(\omega)\}_{\omega \in \Omega}$
 of  $  X$  is called  tempered
 with respect to $(\theta_t)_{t\in \R}$ if for $P$-a.e. $\omega \in \Omega$,
 $$ \lim_{t \to \infty} e^{- \beta t} d(B(\theta_{-t} \omega)) =0
 \quad \mbox{for all} \  \beta>0,
 $$
 where $d(B) =\sup_{x \in B} \| x \|_{X}$.
\end{defn}

\begin{defn}
Let $\mathcal{D}$ be a collection of  random  subsets of $X$.
Then  $\mathcal{D}$ is called inclusion-closed if 
   $D=\{D(\omega)\}_{\omega \in \Omega} \in {\mathcal{D}}$
and  $\tilde{D}=\{\tilde{D}(\omega) \subseteq X:  \omega \in \Omega\} $
with
  $\tilde{D}(\omega) \subseteq D(\omega)$ for all $\omega \in \Omega$ imply
  that  $\tilde{D} \in {\mathcal{D}}$.
  \end{defn}

\begin{defn}
Let $\mathcal{D}$ be a collection of random subsets of $X$ and
$\{K(\omega)\}_{\omega \in \Omega} \in \mathcal{D}$. Then
$\{K(\omega)\}_{\omega \in \Omega} $ is called an absorbing set of
$\Phi$ in $\mathcal{D}$ if for every $B \in \mathcal{D}$ and
$P$-a.e. $\omega \in \Omega$, there exists $t_B(\omega)>0$ such
that
$$
\Phi(t, \theta_{-t} \omega, B(\theta_{-t} \omega)) \subseteq
K(\omega) \quad \mbox{for all} \ t \ge t_B(\omega).
$$
\end{defn}

\begin{defn}
Let $\mathcal{D}$ be a collection of random subsets of $X$. Then
$\Phi$ is said to be  $\mathcal{D}$-pullback asymptotically
compact in $X$ if  for $P$-a.e. $\omega \in \Omega$,
$\{\Phi(t_n, \theta_{-t_n} \omega,
x_n)\}_{n=1}^\infty$ has a convergent  subsequence  in $X$
whenever
  $t_n \to \infty$, and $ x_n\in   B(\theta_{-t_n}\omega)$   with
$\{B(\omega)\}_{\omega \in \Omega} \in \mathcal{D}$.
\end{defn}

\begin{defn}
Let $\mathcal{D}$ be a collection of random subsets of $X$
and $\{\mathcal{A}(\omega)\}_{\omega \in \Omega} \in  \mathcal{D}$.
Then   $\{\mathcal{A}(\omega)\}_{\omega \in \Omega} $
is called a   $\mathcal{D}$-random   attractor
(or $\mathcal{D}$-pullback attractor)  for
  $\Phi$
if the following  conditions are satisfied, for $P$-a.e. $\omega \in \Omega$,

(i) \  $\mathcal{A}(\omega)$ is compact,  and
$\omega \mapsto d(x, \mathcal{A}(\omega))$ is measurable for every
$x \in X$;

(ii) \ $\{\mathcal{A}(\omega)\}_{\omega \in \Omega}$ is invariant, that is,
$$ \Phi(t, \omega, \mathcal{A}(\omega)  )
= \mathcal{A}(\theta_t \omega), \ \  \forall \   t \ge 0;
$$

(iii) \ \ $\{\mathcal{A}(\omega)\}_{\omega \in \Omega}$
attracts  every  set  in $\mathcal{D}$,  that is, for every
 $B = \{B(\omega)\}_{\omega \in \Omega} \in \mathcal{D}$,
$$ \lim_{t \to  \infty} d (\Phi(t, \theta_{-t}\omega, B(\theta_{-t}\omega)), \mathcal{A}(\omega))=0,
$$
where $d$ is the Hausdorff semi-metric given by
$d(Y,Z) =
  \sup_{y \in Y }
\inf_{z\in  Z}  \| y-z\|_{X}
 $ for any $Y\subseteq X$ and $Z \subseteq X$.
\end{defn}

The following existence result on a    random attractor
for a  continuous  RDS
can be found in \cite{bat1,  fla1}.
\begin{prop}
\label{att} Let $\mathcal{D}$ be an inclusion-closed 
 collection of random subsets of
$X$ and $\Phi$ a continuous RDS on $X$ over $(\Omega, \mathcal{F},
P,  (\theta_t)_{t\in \R})$. Suppose  that $\{K(\omega)\}_{\omega
\in K} $ is a closed  absorbing set of  $\Phi$  in $\mathcal{D}$
and $\Phi$ is $\mathcal{D}$-pullback asymptotically compact in
$X$. Then $\Phi$ has a unique $\mathcal{D}$-random attractor
$\{\mathcal{A}(\omega)\}_{\omega \in \Omega}$ which is given by
$$\mathcal{A}(\omega) =  \bigcap_{\tau \ge 0} \  \overline{ \bigcup_{t \ge \tau} \Phi(t, \theta_{-t} \omega, K(\theta_{-t} \omega)) }.
$$
\end{prop}

In this paper, we will 
determine a collection of random subsets
for
the stochastic \bbme on $Q$, and prove
the equation  has a $\mathcal{D}$-random attractor in $\hone$.

\section{Stochastic \bbm equations
  }
\setcounter{equation}{0}

In this section, we discuss the existence of a continuous random dynamical system
for the stochastic \bbme   defined  
 on  an unbounded   channel.
  Let  $D$ be  a
bounded open subset of $\rtwo$   and 
  $Q=D \times \rone$.
 Consider the stochastic   BBM  equation    defined on $Q$:
\be
\label{ue}
d u -d (\De u) - \nu \De u dt + \nabla \cdot {\overrightarrow{F}}(u)  dt
= g dt + h  dw,   \quad x \in Q, \quad t>0,
\ee
with the boundary condition
\be
\label{ue2}
u|_{\partial Q} =0,
\ee
and the initial condition
\be
\label{ue3}
u(x,0) = u_0(x), \quad x \in Q,
\ee
 where $\nu$ is a positive constant, $g \in L^2(Q)$ and
  $h \in H^1_0 (Q)$ are given,  
  $w$   is a    two-sided real-valued
 Wiener process  on a probability
 space which will be specified below,
 and
   $\vecf $   is a  smooth  nonlinear vector function
   given  by
    $\vecf (s)
= (F_{1}(s), F_{2}(s), F_{3}(s))$ for $s \in \rone$, where
  $F_{k} \  (k = 1, 2, 3)$  satisfies
\be
\label{fcond}
F_{k} (0) = 0, \quad  |{ F^\prime_{k}(s)}|  \leq \gamma_{1} + \gamma_{2} |s|, \quad s \in \rone,
\ee
where $\gamma_1$ and $\gamma_2$ are positive constants.
 Denote by
\be
\label{gdef}
 G_{k} (s) = \int_0^s F_{k} (t) dt \quad  \mbox{and}   \quad  
 \vecg (s)  = (G_1(s), G_2(s), G_3(s) ), \quad s \in \rone.
\ee
Then it  follows from \eqref{fcond} that,  for $k = 1, 2, 3,$
\be
\label{fgcond}
 | {F_{k}(s)}|  \leq \gamma_{1}|{s}|  + \gamma_{2} |{s}|^{2}
\quad  \mbox{and} \quad  | G_k (s)|  \leq \gamma_{1} |s|^{2} + \gamma_{2} |s|^{3}.
\ee

Note that   the classical  \bbme with $F_k(s) = s + {\frac 12} s^2$
 indeed satisfies condition \eqref{fcond}.
 Let $\beta_0$ and $\lambda$ be the positive constants
 in \eqref{agmon} and \eqref{poincare}, respectively, and
 denote by
 \be
 \label{fixbeta}
 \delta = \min\{\nu, {\frac 14} \nu \lambda\}
 \quad \mbox{and} \quad 
 \beta = 4 \beta_0 \gamma_2 \| h \|_{H^1}.
 \ee
 Then choose a sufficiently large number $\alpha$
 such that
 \be
 \label{fixalpha}
 \alpha  > {\frac {128 \beta ^2}{\delta ^2}}.
 \ee
 As we will  see later, these constants
 $\alpha, \beta$ and $\delta$  prove useful
when  deriving uniform estimates on the solutions.

In the sequel, we consider the probability space
$(\Omega, \mathcal{F}, P)$ where
$$
\Omega = \{ \omega   \in C(\R, \R ): \ \omega(0) =  0 \},
$$
$\mathcal{F}$ is the Borel $\sigma$-algebra induced by the compact-open topology
of $\Omega$, and $P$ the corresponding Wiener measure on
$(\Omega, \mathcal{F})$.  
 Define the time shift by
$$ \theta_t \omega (\cdot) = \omega (\cdot +t) - \omega (t), \quad  \omega \in \Omega, \ \ t \in \R .
$$
Then $(\Omega, \mathcal{F}, P, (\theta_t)_{t\in \R})$
is a metric dynamical  system. For our purpose, we need
to convert the stochastic equation \eqref{ue} with  a random 
term  into  a deterministic one with  a random parameter. 
To this end, we consider the  
  stationary solutions of the one-dimensional
equation:
\be
\label{y}
dy  + \alpha  y  dt = d w (t),
\ee
where $\alpha$ satisfies \eqref{fixalpha}.
The solution  to \eqref{y}  is   given by
$$
y  (\theta_t \omega )= - \alpha  \int^0_{-\infty} e^{\alpha \tau}  ( \theta_t \omega  )
(\tau) d \tau, \quad t \in \R.
$$
It is known  that  there exists a $\theta_t$-invariant set $\tilde{\Omega}\subseteq \Omega$
of full $P$ measure  such that 
  $y(\theta_t\omega)$  is 
 continuous in $t$ for every $\omega \in \tilde{\Omega}$,
and 
the random variable $|y(\omega)|$ is tempered 
(see, e.g., \cite{arn1, bat1, cra1, cra2}).
   
Put $z (\theta_t \omega) =   (I-\De )^{-1} h   y (\theta_t \omega  ) $
where $\De$ is the Laplacian with domain $H^1_0(Q)\bigcap H^2(Q)$.
By \eqref{y}
we  find that
$$  dz -d (\De z) + \alpha (z- \De z ) dt = h d w.$$
Let $v(t, \omega)  = u(t, \omega) - z(\theta_t \omega)$,
where $u(t, \omega)$  satisfies \eqref{ue}-\eqref{ue3}.
Then for  $v(t, \omega)$ we have that 
\be
\label{ve}
v_t - \De v_t - \nu \De v = -\nabla \cdot \vecf (v+z(\theta_t \omega)) 
+ g + \alpha z(\theta_t \omega) + (\nu -\alpha) \De z(\theta_t \omega),
\ee
with the boundary condition
\be
\label{ve2}
v|_{\partial Q} =0,
\ee
and the initial condition
\be
\label{ve3}
v(0, \omega) = v_0(\omega)  .
\ee

By a Galerkin method as in \cite{med-mil, med-per}, it  can be proved  that under the assumption \eqref{fcond},  
for $P$-a.e. $\omega \in \Omega $ and for all $v_0 \in H^1_0(Q)$,
 problem  \eqref{ve}-\eqref{ve3}
has   a unique solution
$v (\cdot, \omega, v_0) 
  \in C([0, \infty), \hone)  $
  with $v(0, \omega, v_0) = v_0$.   
  Further,  the solution $v (t, \omega, v_0) $
  is continuous    with respect to $v_0 $ in $\hone$
for all $t \ge 0$.  Throughout this paper, we always write
\be
\label{uv}
 u(t, \omega, u_0)= 
    v (t, \omega, v_0) + z (\theta_t \omega),
\quad \mbox{with} \quad   v_0 =  u_0 -
z (\omega) .
\ee 
Then $u$ is  a
solution of problem \eqref{ue}-\eqref{ue3} in some sense. We now
define  a mapping  $\Phi: \R^+ \times \Omega \times \hone  \to \hone$ by 
\be \label{phi}
 \Phi (t, \omega, u_0  ) =   u (t, \omega, u_0) , 
 \quad \forall \ (t, \omega, u_0) \in
\R^+ \times \Omega \times\hone. \ee 
Note that
$\Phi$   satisfies conditions (i), (ii) and (iii) in Definition
\ref{RDS}. Therefore, $\Phi$ is a continuous  random dynamical
system
 associated  with  the stochastic  \bbme
    on $Q$.  In what   follows, we will prove 
that $\Phi$ has a  $\mathcal{D}$-random attractor 
in $\hone$, where 
$\mathcal{D}$ is a collection of random subsets of 
$\hone$ given by
\be
\label{Dset}
 \mathcal{D} = \{ B: \ B= \{B(\omega)\}_{\omega\in \Omega},  \   B(\omega) \subseteq H^1_0(Q)
 \ \mbox{and } \
 \ e^{-{\frac 18}\delta   t} d(B(\theta_{-t}\omega)) \to 0 \ \mbox{as} \ t \to \infty \},
\ee
 where $\delta$  is the  positive constant 
 in \eqref{fixbeta}  and
$$   d(B(\theta_{-t}\omega)) =\sup_{u \in B(\theta_{-t}\omega)} \| u \|_{\hone}.$$
 Notice that $\mathcal{D}$ contains all tempered random sets,
 especially all bounded  deterministic subsets of $H^1_0(Q)$.

\section{Uniform  estimates   }
\setcounter{equation}{0}

In this section, we
 derive uniform estimates on the  solutions of the  stochastic
\bbme  defined on $Q$
when $t \to \infty$, which include the uniform estimates on the
tails of solutions  as both $x$ and $t$ approach infinity.
   These estimates  are necessary
 for proving  the existence of bounded absorbing sets
 and the asymptotic compactness of
  the random dynamical system.

 From now on, we always assume that $\mathcal{D}$ is
 the collection of random subsets of $\hone$ given
 by \eqref{Dset}. We first  derive the following uniform estimates
 on $v$ in $\hone$.

\begin{lem}
\label{lem41}
 Assume that $g \in \h$, $ h  \in  H^1_0(Q)$ and
 \eqref{fcond} holds. Let
 $B=\{B(\omega)\}_{\omega \in \Omega}\in \mathcal{D}$ and
 $v_0(\omega) \in B(\omega)$. Then for $P$-a.e. $\omega \in \Omega$,
 there is $T= T (B, \omega) >0$ such that for all $t \ge T $,
 $$
 \| v(t, \theta_{-t} \omega, v_0(\theta_{-t} \omega)  )\|_{\hone}
 \le r_1(\omega),
 $$
 where $r_1(\omega)$ is a positive random function satisfying
\be
\label{lem41_1}
 e^{- {\frac 18}\delta t} r_1(\theta_{-t} \omega) \to 0
 \quad \mbox{as} \quad t \to \infty.
\ee
 \end{lem}
 
 \begin{proof}
 Taking the inner product of \eqref{ve} with $v$   in  $\h$ we find that
$$
 {\frac 12} {\frac d{dt}} \left ( \| v \|^2 + \| \nabla v \|^2 
 \right ) + \nu \| \nabla v \|^2
 = - \int_Q v \nabla \cdot \vecf (v+\zt)  dx
 $$
  \be
 \label{p41_1}
 + (g + \alpha \zt + (\nu -\alpha) \De \zt, v).
 \ee
 By \eqref{gdef} we have 
 $\nabla \cdot \vecg(u) = \vecf (u) \cdot \nabla u$ and hence,
 by \eqref{agmon} and \eqref{fgcond}, 
 the nonlinear term on the right-hand side of \eqref{p41_1}
 satisfies
 $$-\int_Q v \nabla \cdot \vecf(v +\zt) dx
 =-\int_Q (u-\zt) \nabla \cdot \vecf(v +\zt) dx
 $$
 $$=
 -\int_Q u \nabla \cdot \vecf(u) dx
 +  \int_Q \zt  \nabla \cdot \vecf(v +\zt) dx
 $$
 $$=
  \int_Q   \vecf(u) \cdot \nabla u dx
 +  \int_Q \zt  \nabla \cdot \vecf(v +\zt) dx
 $$
 $$=
  \int_Q     \nabla \cdot  \vecg (u) dx
 +  \int_Q \zt  \nabla \cdot \vecf(v +\zt) dx
 $$
 $$
 = \int_Q \zt  \nabla \cdot \vecf(v +\zt) dx
 =- \int_Q     \vecf(v +\zt) \cdot \nabla \zt dx
 $$
 $$
 \le \gamma_1 \int_Q   |v + \zt| \ |\nabla \zt | dx
 + \gamma_2 \int_Q | v + \zt|^2   |\nabla \zt | dx
 $$
 $$
 \le
 \gamma_1(  \|v \|  
 +   \| \zt \| )  \| \nabla \zt \|
 + 2\gamma_2  \int_Q |v|^2 |\nabla \zt | dx
 + 2\gamma_2 \int_Q |\zt|^2 |\nabla \zt|dx
 $$
 $$\le {\frac 18}\nu \lambda \| v \|^2
 +c_1 |y(\theta_t \omega)|^2
 + 2\gamma_2 \| \nabla \zt \|_\infty \| v\|^2
 + 2\gamma_2 \| \zt \|^2_4 \| \nabla  \zt \|
 $$
 $$\le {\frac 18}\nu \lambda \| v \|^2
 +c_1 |y(\theta_t \omega)|^2
 + 2\gamma_2 \beta_0  \| \nabla \zt \|_{H^2} \| v\|^2
 + c_2 | y(\theta_t \omega) |^3
 $$
 \be  
\label{p41_2}
\le {\frac 18}\nu \lambda \| v \|^2
 +c_1 |y(\theta_t \omega)|^2
 + 2\gamma_2 \beta_0 \|h\|_{H^1}   | y(\theta_t \omega) | \ \| v\|^2
 + c_2 | y(\theta_t \omega) |^3.
\ee
By the Young inequality, the second term on the right-hand side 
of \eqref{p41_1} is bounded by
$$
|(g + \alpha \zt + (\nu -\alpha) \De \zt, v)|
\le {\frac 18} \nu \lambda \| v \|^2
+ c_3 (\| g \|^2 + \| \zt \|^2 +  \|  \De \zt \|^2 )
$$
\be
\label{p41_3}
\le {\frac 18} \nu \lambda \| v \|^2 
+ c_4 (1 + |y(\theta_t  \omega )| ^2 ).
\ee
It follows from \eqref{p41_1}-\eqref{p41_3} that
\be
\label{p41_4}
{\frac d{dt}} \| v \|^2_{H^1 }
+ 2 \nu \| \nabla v \|^2
\le  {\frac 12} \nu \lambda \| v \|^2 + 
4 \beta_0 \gamma_2 \| h \|_{H^1} |y(\theta_t \omega)|  \ \| v \|^2
+ c (1 + |y(\theta_t \omega)|^2 + |y(\theta_t \omega)|^3 ).
\ee
By \eqref{poincare} we have that
\be
\label{p41_5}
2 \nu \| \nabla v \|^2 \ge \nu \| \nabla v \|^2 + \nu \lambda \| v \|^2.
\ee
By \eqref{p41_4} and \eqref{p41_5} we get 
$$
{\frac d{dt}} \| v \|^2_{H^1 }
+  \nu \| \nabla v \|^2
+  {\frac 12} \nu \lambda \| v \|^2 \le 
4 \beta_0 \gamma_2 \| h \|_{H^1} |y(\theta_t \omega)|  \ \| v \|^2
+ c (1 + |y(\theta_t \omega)|^2 + |y(\theta_t \omega)|^3 ),
$$
which along with \eqref{fixbeta} implies that
\be
\label{p41_6}
{\frac d{dt}} \| v \|^2_{H^1 }
+ ( \delta  - \beta  |y(\theta_t \omega)| ) \| v \|^2_{H^1}
 \le   
  c (1 + |y(\theta_t \omega)|^2 + |y(\theta_t \omega)|^3 ).
\ee
Multiplying \eqref{p41_6} by $e^{\int_0^t(\delta -\beta |y(\theta_\tau \omega)| ) d\tau}$
and then integrating over $(0,s)$ with $s \ge 0$, we obtain that
$$
\| v(s, \omega, v_0(\omega ) )\|^2_{H^1}
\le e^{ -\delta s + \beta \int_0^s  |y(\theta_\tau \omega)|  d \tau} \| v_0 (\omega) \|^2_{H^1}
$$
\be\label{p41_7}
+ c  \int_0^s (1 + |y(\theta_\sigma \omega)|^2 + |y(\theta_\sigma \omega)|^3 ) 
e^{  \delta (\sigma -s) + \beta \int_\sigma^s  |y(\theta_\tau \omega)|  d \tau} d \sigma.
\ee
We now replace $\omega$ by $\theta_{-t} \omega$ with $t \ge 0$ in \eqref{p41_7} to
get that, for any $s \ge 0$ and $t \ge 0$, 
$$
\| v(s, \theta_{-t}\omega, v_0(\theta_{-t}\omega ) )\|^2_{H^1}
\le e^{ -\delta s + \beta \int_0^s  |y(\theta_{\tau-t} \omega)|  d \tau} \| v_0 (\theta_{-t}\omega) \|^2_{H^1}
$$
$$
+ c  \int_0^s (1 + |y(\theta_{\sigma-t} \omega)|^2 + |y(\theta_{\sigma-t} \omega)|^3 ) 
e^{  \delta (\sigma -s) + \beta \int_\sigma^s  |y(\theta_{\tau-t} \omega)|  d \tau} d \sigma.
$$
$$
\le e^{ -\delta s + \beta \int_{-t}^{s-t}  |y(\theta_{\tau} \omega)|  d \tau} \| v_0 (\theta_{-t}\omega) \|^2_{H^1}
$$
 \be
\label{vs}
+ c  \int_{-t}^{s-t} (1 + |y(\theta_{\sigma} \omega)|^2 + |y(\theta_{\sigma} \omega)|^3 ) 
e^{  \delta (\sigma -s +t) + \beta \int_\sigma^{s-t}  |y(\theta_{\tau} \omega)|  d \tau} d \sigma.
\ee
By \eqref{vs} we find that, for all $t\ge 0$,
$$
\| v(t, \theta_{-t}\omega, v_0(\theta_{-t}\omega ) )\|^2_{H^1}
\le e^{ -\delta t + \beta \int_{-t}^{0}  |y(\theta_{\tau} \omega)|  d \tau} \| v_0 (\theta_{-t}\omega) \|^2_{H^1}
$$
\be
\label{p41_8}
+ c  \int_{-t}^{0} (1 + |y(\theta_{\sigma} \omega)|^2 + |y(\theta_{\sigma} \omega)|^3 ) 
e^{  \delta \sigma  + \beta \int_\sigma^{0}  |y(\theta_{\tau} \omega)|  d \tau} d \sigma.
\ee
Note that $|y(\theta_\tau \omega)|$ is stationary and ergodic
(see, e.g. \cite{cra1}).  Then it follows from the ergodic theorem
that
$$
\lim_{t \to \infty} {\frac 1t} \int_{-t}^0 | y(\theta_\tau \omega) | d\tau
= E(|y(\omega)|).
$$
On the other hand, we have
$$
E(|y(\omega)|) \le \left ( E(|y(\omega)|^2) \right )^{\frac 12}
\le {\frac 1{\sqrt{2 \alpha}}},
$$
which shows that
\be
\label{p41_9}
\lim_{t \to \infty} {\frac 1t} \int_{-t}^0 | y(\theta_\tau \omega) | d\tau
< {\frac 2{\sqrt{2 \alpha}}}.
\ee
By \eqref{fixalpha} and \eqref{p41_9} we find that there is $T_0(\omega)>0$
such that for all $t \ge T_0(\omega)$,
\be
\label{deltacond}
 \beta \int_{-t}^{0}  |y(\theta_{\tau} \omega)|  d \tau
 <  {\frac {2\beta t}{\sqrt{2\alpha}}}  < {\frac 18} \delta t.
 \ee
 By \eqref{p41_8} and \eqref{deltacond} we find that, for all
 $t\ge T_0(\omega)$,
 $$
\| v(t, \theta_{-t}\omega, v_0(\theta_{-t}\omega ) )\|^2_{H^1}
\le e^{ - {\frac 78}\delta t }  \| v_0 (\theta_{-t}\omega) \|^2_{H^1}
$$
\be
\label{p41_10}
+ c  \int_{-t}^{0} (1 + |y(\theta_{\sigma} \omega)|^2 + |y(\theta_{\sigma} \omega)|^3 ) 
e^{  \delta \sigma  + \beta \int_\sigma^{0}  |y(\theta_{\tau} \omega)|  d \tau} d \sigma.
\ee
Note that $|y(\theta_{\sigma} \omega)|$ is tempered, and hence by
\eqref{deltacond},  the integrand of the second term on the right-hand side of
\eqref{p41_10} is convergent to zero exponentially  as $ \sigma \to  - \infty$. 
This shows that the following integral is convergent:
\be
\label{p41_11}
r_0(\omega) = c  \int_{- \infty}^{0} (1 + |y(\theta_{\sigma} \omega)|^2 + |y(\theta_{\sigma} \omega)|^3 ) 
e^{  \delta \sigma  + \beta \int_\sigma^{0}  |y(\theta_{\tau} \omega)|  d \tau} d \sigma.
\ee
It follows from \eqref{p41_10}-\eqref{p41_11} that, for all $t \ge T_0(\omega)$,
\be
\label{p41_20}
\| v(t, \theta_{-t}\omega, v_0(\theta_{-t}\omega ) )\|^2_{H^1}
\le e^{ - {\frac 78}\delta t }  \| v_0 (\theta_{-t}\omega) \|^2_{H^1}
+ r_0(\omega).
\ee
On the other hand, by assumption, $B=\{B(\omega)\}_{\omega \in \Omega} \in \mathcal{D}$
and hence  we have
$$ e^{ - {\frac 14}\delta t }  \| v_0 (\theta_{-t}\omega) \|_{H^1} \to 0
\quad \mbox{as} 
\quad t \to \infty,$$
from which and \eqref{p41_20} we find that  there is
$T=T(B,\omega)>0$ such that for all $t\ge T $,
$$
\| v(t, \theta_{-t}\omega, v_0(\theta_{-t}\omega ) )\|^2_{H^1}
\le   2 r_0(\omega).
$$
Let $r_1(\omega) =\sqrt{2r_0(\omega)}$. Then we get that, for all
$t\ge T $,
\be
\label{p41_21}
\| v(t, \theta_{-t}\omega, v_0(\theta_{-t}\omega ) )\| _{H^1}
\le    r_1(\omega).
\ee
Next, we prove $r_1(\omega)$ satisfies \eqref{lem41_1}.
Replacing $\omega$ by $\theta_{-t} \omega$ in \eqref{p41_11} we obtain that
$$
r_0(\theta_{-t}\omega) = 
c  \int_{- \infty}^{0} (1 + |y(\theta_{\sigma -t} \omega)|^2 + |y(\theta_{\sigma-t} \omega)|^3 ) 
e^{  \delta \sigma  + \beta \int_\sigma^{0}  |y(\theta_{\tau -t} \omega)|  d \tau} d \sigma
$$
 $$  = 
c  \int_{- \infty}^{-t} (1 + |y(\theta_{\sigma } \omega)|^2 + |y(\theta_{\sigma} \omega)|^3 ) 
e^{  \delta (\sigma+t)  + \beta \int_\sigma^{-t}  |y(\theta_{\tau } \omega)|  d \tau} d \sigma.
$$
 $$ \le
c  \int_{- \infty}^{-t} (1 + |y(\theta_{\sigma } \omega)|^2 + |y(\theta_{\sigma} \omega)|^3 ) 
e^{ {\frac 3{16}} \delta (\sigma+t)  + \beta \int_\sigma^{-t}  |y(\theta_{\tau } \omega)|  d \tau} d \sigma.
$$
 $$ \le
c  e^{{\frac 3{16}}\delta t} \int_{- \infty}^{-t} (1 + |y(\theta_{\sigma } \omega)|^2 + |y(\theta_{\sigma} \omega)|^3 ) 
e^{ {\frac 3{16}} \delta \sigma  + \beta \int_\sigma^{0}  |y(\theta_{\tau } \omega)|  d \tau} d \sigma.
$$
 \be
\label{p41_30}
 \le
c  e^{{\frac 3{16}}\delta t} \int_{- \infty}^{0} (1 + |y(\theta_{\sigma } \omega)|^2 + |y(\theta_{\sigma} \omega)|^3 ) 
e^{ {\frac 3{16}} \delta \sigma  + \beta \int_\sigma^{0}  |y(\theta_{\tau } \omega)|  d \tau} d \sigma.
\ee
Note that the last integral in the above is indeed convergent since
the integrand converges to zero exponentially by \eqref{deltacond}.
Then we have
$$
e^{-{\frac 18} \delta t} r_1(\theta_{-t} \omega)
= 
e^{-{\frac 18} \delta t}  \sqrt{2 r_0(\theta_{-t} \omega)}
$$
$$\le
\sqrt{2c} e^{-{\frac 1{32}} \delta t} 
 \left (  \int_{- \infty}^{0} (1 + |y(\theta_{\sigma } \omega)|^2 + |y(\theta_{\sigma} \omega)|^3 ) 
e^{ {\frac 3{16}} \delta \sigma  + \beta \int_\sigma^{0}  |y(\theta_{\tau } \omega)|  d \tau} d \sigma
\right )^{\frac 12} \to 0, \quad \mbox{as} \ t \to \infty,
$$
which along with \eqref{p41_21} completes the proof.
 \end{proof}

\begin{lem}
\label{lem42}
 Assume that $g \in \h$, $ h  \in  H^1 _0 (Q)$ and
 \eqref{fcond} holds. Let
 $B=\{B(\omega)\}_{\omega \in \Omega}\in \mathcal{D}$ and
 $v_0(\omega) \in B(\omega)$. Then for $P$-a.e. $\omega \in \Omega$,
   every $s\ge 0$ and $t\ge 0$, we have 
 $$
\| v_s (s, \theta_{-t}\omega, v_0(\theta_{-t}\omega ) )\|^2_{H^1}
\le c+  c e^{ - 2 \delta s + 2 \beta \int_{-t}^{s-t}  |y(\theta_{\tau} \omega)|  d \tau}
 \| v_0 (\theta_{-t}\omega) \|^4_{H^1}
$$
$$+ c \left ( \int_{-t}^{s-t} (1 + |y(\theta_{\sigma} \omega)|^2 + |y(\theta_{\sigma} \omega)|^3 ) 
e^{  \delta (\sigma -s +t) + \beta \int_\sigma^{s-t}  |y(\theta_{\tau} \omega)|  d \tau} d \sigma
\right )^2
$$
$$
+
c \left (
\| z(\theta_{s-t} \omega)\|^2
+ \| z(\theta_{s-t} \omega)\|^4_{H^1}
+\| z(\theta_{s-t} \omega)\|^2_{H^2}
\right ),
$$
where $c$ is a positive deterministic constant.
 \end{lem}
 
 \begin{proof} Taking the inner product of \eqref{ve}
 with $v_t$ in $\h$ we obtain that
 \be
 \label{p42_1}
 \| v_t  \| ^2  + \| \nabla v_t \|^2
 +  \nu (\nabla v, \nabla v_t)
 = \int_Q     \vecf (v+z(\theta_t \omega ) )\cdot \nabla v_t dx
 +
( g + \alpha z(\theta_t \omega) + (\nu -\alpha) \De z(\theta_t \omega), v_t ).
\ee
We now estimate every term in the above. First we have
\be
\label{p42_2}
\nu |(\nabla v, \nabla v_t)|
\le \nu \| \nabla v \| \ \| \nabla v_t \| 
\le {\frac 14}\| \nabla v_t \|^2 + \nu^2 \| \nabla v \|^2.
\ee
By \eqref{fgcond}, the nonlinear term in \eqref{p42_1} is bounded by
$$
| \int_Q     \vecf (v+z(\theta_t \omega ) )\cdot \nabla v_t dx|
\le \gamma_1 \int_Q |v+z(\theta_t \omega ) | \ |\nabla v_t | dx
+ \gamma_2 \int_Q |v+z(\theta_t \omega ) | ^2 \ |\nabla v_t | dx
$$
$$
\le {\frac 14} \| \nabla v_t \|^2
+c (\| v\|^2 + \| v \|_4^4 + \| z(\theta_t \omega )\|^2
+ \| z(\theta_t \omega )\|^4_4 )
$$
\be
\label{p42_3}
\le {\frac 14} \| \nabla v_t \|^2
+c (\| v\|^2 + \| v \|_{H^1}^4 + \| z(\theta_t \omega )\|^2
+ \| z(\theta_t \omega )\|^4_{H^1} ).
\ee
For the last term on the right-hand side of \eqref{p42_1} we have
\be
\label{p42_4}
| ( g + \alpha z(\theta_t \omega) + (\nu -\alpha) \De z(\theta_t \omega), v_t )|
\le {\frac 12} \| v_t \|^2
+ c( \| g \|^2 + \|z(\theta_t \omega )\|^2 +\|\De z(\theta_t \omega )\|^2).
\ee
Then it follows from \eqref{p42_1}-\eqref{p42_4} that
$$
\|v_t \|^2 + \| \nabla v_t \|^2
\le  
  c_1  (\| v \|^2_{H^1} + \| v \|_{H^1}^4)
+ c_1 (1 + \|z(\theta_t \omega )\|^2
+ \|z(\theta_t \omega )\|^4_{H^1} + 
\|z(\theta_t \omega )\|^2_{H^2} )
$$
$$
\le  
  c   \| v \|_{H^1}^4 
+ c(1 + \|z(\theta_t \omega )\|^2
+ \|z(\theta_t \omega )\|^4_{H^1} + 
\|z(\theta_t \omega )\|^2_{H^2} ),
$$
which shows that, for all $t \ge 0$,
\be
\label{p42_5}
\| v_t (t, \omega, v_0(\omega ) )\|^2_{H^1}
\le  
  c   \| v  (t, \omega, v_0(\omega ) )  \|_{H^1}^4 
+ c(1 + \|z(\theta_t \omega )\|^2
+ \|z(\theta_t \omega )\|^4_{H^1} + 
\|z(\theta_t \omega )\|^2_{H^2} ).
\ee
First replacing $t$ by $s$ and then replacing 
$\omega$ by $\theta_{t} \omega$ in \eqref{p42_5},
 we get that,
for all $s \ge 0$ and $t \ge 0$,
$$
\| v_s (s, \theta_{-t}\omega, v_0(\theta_{-t}\omega ) )\|^2_{H^1}
$$
$$
\le  
  c   \| v  (s, \theta_{-t}\omega, v_0(\theta_{-t}\omega ) )  \|_{H^1}^4 
+ c(1 + \|z(\theta_{s-t} \omega )\|^2
+ \|z(\theta_{s-t} \omega )\|^4_{H^1} + 
\|z(\theta_{s-t} \omega )\|^2_{H^2} ),
$$
which along with \eqref{vs} completes the proof.
 \end{proof}
 
We are now ready to derive the uniform estimates
on the tails of solutions when $x$ and $t$ approach infinity,
which are crucial for proving the asymptotic compactness
of the equation.
To this end,  for every $x \in Q = \Omega\times\rone$, we  will
write
$x = (x_{1},x_{2},x_{3})$ where $(x_{1},x_{2}) \in \Omega$ and $x_{3} \in \rone$.
Given $k>0$, denote by $Q_{k} =$\{$(x_{1},x_{2},x_{3}) \in Q$: $|{x_{3}}|<k\}$,
 and $Q \backslash Q_{k}$   the complement of $Q_{k}$.

\begin{lem}
\label{lem43}
 Assume that $g \in \h$, $ h  \in  H^1_0 (Q)$ and
 \eqref{fcond} holds. Let
 $B=\{B(\omega)\}_{\omega \in \Omega}\in \mathcal{D}$ and
 $v_0(\omega) \in B(\omega)$. Then for 
every $\epsilon>0$ and $P$-a.e. $\omega \in \Omega$,
  there exist $T=T(B, \omega, \epsilon)>0$ and $k_0=k_0(\omega, \epsilon)>0$
  such that for all $t \ge T$,
  $$
  \int_{Q\backslash Q_{k_0}} \left (
   |v(t, \theta_{-t} \omega, v_0(\theta_{-t} \omega ) )|^2
   +  |\nabla v(t, \theta_{-t} \omega, v_0(\theta_{-t} \omega ) )|^2
  \right ) dx \le \epsilon.
  $$
 \end{lem}
 
 \begin{proof}
  Take a smooth function $\phi$ such that $0 \leq \phi \leq 1$   for all  $s \in \rone$ and
\be
\label{phif}
\phi (s)  =
\left \{
\begin{array}{ll}
    0, & \quad  \mbox{if }  \quad   |s| < 1,\\
  1, & \quad  \mbox{if }   \quad |s| > 2.
\end{array}
\right.
\ee
Then there is a positive constant $c$ such that
$|\phi^\prime (s)| + |{\phi^{\prime\prime}}  (s) | \le c$
for all $s\in \R$.
Multiplying \eqref{ve}  by $\phi^{2}({\frac{x_{3}^{2}}{k^{2}}})v$ and then integrating
with respect to $x$ on $Q$, we get
$$
\int_Q \cutf vv_t  \ dx - \int_Q \cutf v \De v_t \   dx
- \nu  \int_Q \cutf v \De v \   dx
$$
\be
\label{p43_1}
= - \int_Q \cutf v \nabla \cdot \vecf (v+z(\theta_t \omega))  \ dx
+ \int_Q \cutf v  ( g + \alpha z(\theta_t \omega) + (\nu -\alpha) \De z(\theta_t \omega) )  dx.
\ee
We now deal with the left-hand side of the above.
For the first term on the left-hand side of \eqref{p43_1} we have
\be 
\label{p43_2}
\int_Q \cutf vv_t  \ dx  = {\frac 12} {\frac d{dt}} \int_Q \cutf  |v|^2 dx.
\ee
We also have
$$
- \int_Q \cutf v \De v_t \   dx
= \int_Q \cutf (\nabla v_t \cdot \nabla v) dx
+ \int_Q v \left (\nabla v_t \cdot \nabla   \cutf \right ) dx
$$
\be
\label{p43_3} =
{\frac 12} {\frac d{dt}} \int_Q\cutf |\nabla v|^2 dx
+ \int_Q v \left (\nabla v_t \cdot \nabla   \cutf \right ) dx.
\ee
The last term on the left-hand side of \eqref{p43_1} satisfies
\be
\label{p43_4}
- \nu  \int_Q \cutf v \De v \   dx
  =\nu \int_Q \cutf |\nabla v|^2 dx
+ \nu \int_Q v \left (\nabla v \cdot \nabla \cutf \right ) dx.
\ee
Then it follows from \eqref{p43_1}-\eqref{p43_4} that
$$
{\frac 12} {\frac d{dt}} \int_Q \cutf (|v|^2 + |\nabla v|^2 ) dx
+ \nu \int_Q \cutf |\nabla v |^2 dx
$$
$$
= -\int_Q v \left ( \nabla v_t \cdot \nabla \cutf \right ) dx
-\nu \int_Q v \left (
\nabla v \cdot \nabla \cutf
\right ) dx
$$
$$
+\int_Q\cutf \left ( \vecf (v+z(\theta_t \omega))  \cdot \nabla v \right ) dx
+ \int_Q v \left ( \vecf (v+z(\theta_t \omega)) \cdot \nabla \cutf
\right ) dx
$$
\be
\label{p43_10}
+ \int_Q \cutf v  ( g + \alpha z(\theta_t \omega) + (\nu -\alpha) \De z(\theta_t \omega) )  dx.
\ee
Next, we estimate the right-hand side of \eqref{p43_10}.
For the first term we have
$$
| \int_Q v \left ( \nabla v_t \cdot \nabla \cutf \right ) dx|
\le \int_Q |v| \ |\nabla v_t| \ |2  \phi   \phi^\prime \left ({\frac {x_3^2}{k^2}}
\right ) | \ {\frac {2 |x_3|}{k^2}} dx
$$
$$
\le \int_{k \le |x_3| \le {\sqrt{2}} k} |v| \ |\nabla v_t| \ |2  \phi   \phi^\prime \left ({\frac {x_3^2}{k^2}}
\right ) | \ {\frac {2 |x_3|}{k^2}} dx
\le {\frac ck}  \int_{k \le |x_3| \le {\sqrt{2}} k}  |v| \ |\nabla v_t|   dx
$$
\be
\label{p43_11}
\le {\frac ck} \| v\| \ \| \nabla v_t \|
\le {\frac ck} \| \nabla v_t \|^2 + {\frac ck} \| v \|^2.
\ee
Similarly, the second term on the right-hand side of \eqref{p43_10}
is bounded by
\be
\label{p43_12}
\nu | \int_Q v \left (
\nabla v \cdot \nabla \cutf
\right ) dx|
\le {\frac ck} \| \nabla v  \|^2 + {\frac ck} \| v \|^2.
\ee
For the third term on the right-hand side of \eqref{p43_10} we have
$$
\int_Q\cutf \left ( \vecf (v+z(\theta_t \omega))  \cdot \nabla v \right ) dx
$$
$$
=\int_Q\cutf \left ( \vecf (u )  \cdot \nabla u \right ) dx
- \int_Q\cutf \left ( \vecf (v+z(\theta_t \omega))  \cdot \nabla \zt \right ) dx
$$
 $$
=\int_Q\cutf ( \nabla \cdot \vecg (u) )   dx
- \int_Q\cutf \left ( \vecf (v+z(\theta_t \omega))  \cdot \nabla \zt \right ) dx
$$
\be
\label{p43_13}
= - \int_Q \vecg(u) \cdot \nabla \cutf   dx
- \int_Q\cutf \left ( \vecf (v+z(\theta_t \omega))  \cdot \nabla \zt \right ) dx.
\ee
By \eqref{fgcond},  the first term of the above is bounded by
$$
| \int_Q \vecg(u) \cdot \nabla \cutf dx | 
\le \int_Q (\gamma_1 |u|^2 + \gamma_2 |u|^3 )
|2 \phi \phi^{\prime} \left ({\frac {x_3^2}{k^2}}
\right )| \  {\frac {2|x_3|}{k^2}} dx
$$
$$
\le \int_{k\le |x_3| \le \sqrt{2} k} (\gamma_1 |u|^2 + \gamma_2 |u|^3 )
|2 \phi \phi^{\prime} \left ({\frac {x_3^2}{k^2}}
\right )| \  {\frac {2|x_3|}{k^2}} dx
$$
$$
\le {\frac ck} \int_{k\le |x_3| \le \sqrt{2} k} ( \gamma_1   |u|^2 + \gamma_2  |u|^3 )
  dx\le {\frac ck} \| u\|^2 + {\frac ck} \| u \|_3^3
$$
\be
\label{p43_14}
\le {\frac ck} \| v\|^2 + {\frac ck} \| v \|_3^3
+   {\frac ck} \| \zt \|^2 + {\frac ck} \| \zt \|_3^3.
\ee
By \eqref{agmon} and \eqref{fgcond},  for the second term of \eqref{p43_13} we have
$$ | \int_Q\cutf \left ( \vecf (v+z(\theta_t \omega))  \cdot \nabla \zt \right ) dx|
 $$
 $$
 \le \gamma_1 \int_Q  \cutf   |v + \zt| \ |\nabla \zt | dx
 + \gamma_2 \int_Q  \cutf | v + \zt|^2   |\nabla \zt | dx
 $$
 $$
 \le \gamma_1  \int_Q  \cutf (  |v| \ |\nabla \zt |  + | \zt| \ | \nabla \zt | ) dx
 $$
 $$
 + 2 \gamma_2  \int_Q \cutf  ( |v|^2 |\nabla \zt |  + | \zt| ^2 | \nabla \zt | ) dx
 $$
 $$\le {\frac 1{16}}\nu \lambda \int_Q \cutf |v|^2 dx
 + c \int_Q\cutf \left ( |\zt|^2 + |\nabla \zt |^2 \right )dx
 $$
 $$
 + 2 \gamma_2  \| \nabla \zt \|_{\infty} \int_Q \cutf |v|^2 dx
 + c \int_Q\cutf | \zt| ^2 | \nabla \zt | dx
 $$
 $$\le {\frac 1{16}}\nu \lambda \int_Q \cutf |v|^2 dx
 + c \int_Q\cutf \left ( |\zt|^2 + |\nabla \zt |^2 \right )dx
 $$
 $$
 + 2 \gamma_2 \beta_0  \| \nabla \zt \|_{H^2} \int_Q \cutf |v|^2 dx
 + c \int_Q\cutf | \zt| ^2 | \nabla \zt | dx
 $$
 $$\le {\frac 1{16}}\nu \lambda \int_Q \cutf |v|^2 dx
 + c \int_Q\cutf \left ( |\zt|^2 + |\nabla \zt |^2 \right )dx
 $$
\be
\label{p43_15}
 + 2 \gamma_2 \beta_0 \| h \|_{H^1}    | y(\theta_t \omega)  |  \int_Q \cutf |v|^2 dx
 + c \int_Q\cutf | \zt| ^2 | \nabla \zt | dx.
\ee
It follows from \eqref{p43_13}-\eqref{p43_15} that
$$
 | \int_Q\cutf \left ( \vecf (v+z(\theta_t \omega))  \cdot \nabla v \right )dx | 
$$
$$
 \le {\frac ck} \| v\|^2 + {\frac ck} \| v \|_3^3
+   {\frac ck} \| \zt \|^2 + {\frac ck} \| \zt \|_3^3
$$
$$+  {\frac 1{16}}\nu \lambda \int_Q \cutf |v|^2 dx
 + c \int_Q\cutf \left ( |\zt|^2 + |\nabla \zt |^2 \right )dx
 $$
\be
\label{p43_16}
 + 2 \gamma_2 \beta_0 \| h \|_{H^1}    | y(\theta_t \omega)  |  \int_Q \cutf |v|^2 dx
 + c \int_Q\cutf | \zt| ^2 | \nabla \zt | dx.
\ee
By \eqref{fgcond}, the fourth term on the right-hand side of \eqref{p43_10} is bounded
by
$$
|\int_Q v \left ( \vecf (v+z(\theta_t \omega)) \cdot \nabla \cutf
\right ) dx |
$$
$$
\le
\int_Q |\vecf(v+z(\theta_t \omega)) | \  |2 \phi \phi^{\prime} \left (
{\frac {x_3^2}{k^2}} \right ) | \ {\frac {2|x_3|}{k^2}} |v| dx
$$
$$
\le
\int_{k \le |x_3| \le \sqrt{2} k}
 |\vecf(v+z(\theta_t \omega)) | \  |2 \phi \phi^{\prime} \left (
{\frac {x_3^2}{k^2}} \right ) | \ {\frac {2|x_3|}{k^2}} |v| dx
$$
$$
\le {\frac ck} \int_Q (|v + \zt| + |v + \zt |^2 ) |v| dx
$$
\be
\label{p43_17} \le 
{\frac ck} \left (
\| v\|^2 + \| v \|^3_3 + \| \zt \|^2 + \| \zt \|^3_3
\right ).
\ee
By the Young inequality, the last term on the right-hand side 
of \eqref{p43_10} is bounded by
 $$
| \int_Q \cutf v  ( g + \alpha z(\theta_t \omega) + (\nu -\alpha) \De z(\theta_t \omega) )  dx|
$$
\be
\label{p43_18}
 \le  {\frac 1{16}} \nu \lambda   \int_Q \cutf |v|^2 dx
 + c \int_Q \cutf (g^2 + |\zt |^2 + | \De \zt |^2) dx.
\ee
Finally, by \eqref{p43_10}-\eqref{p43_12} and \eqref{p43_16}-\eqref{p43_18}, we find that
$$
  {\frac d{dt}} \int_Q \cutf (|v|^2 + |\nabla v|^2 ) dx
+ 2  \nu \int_Q \cutf |\nabla v |^2 dx
$$
$$
\le \left ( {\frac 14} \nu \lambda + \beta |y(\theta_t \omega )| \right )
\int_Q\cutf |v|^2 dx
+ {\frac ck} (\| v \|^2 + \| \nabla v \|^2 ) + {\frac ck} \| v\|^3_3 + {\frac ck} \| \nabla v_t \|^2
$$
$$
+ c\int_Q \cutf \left (
 g^2 + |\zt|^2 + |\De \zt|^2 + |\nabla \zt |^2 + |\zt |^2 \ |\nabla \zt |
\right ) dx
$$
\be
\label{p43_20}
+ {\frac ck} \left ( \| \zt \|^2 + \| \zt \|^3_3 \right ).
\ee
We now deal with the second term on the left-hand side of the above. Note that
$$
\int_{Q}|{\nabla \left (\phi({\frac{x_{3}^{2}}{k^{2}}})v \right )}|^{2}dx  =
 \int_{Q}| {v \nabla\phi ({\frac{x_{3}^{2}}{k^{2}}}) + \phi({\frac{x_{3}^{2}}{k^{2}}}) \nabla v}|
^{2}dx  
$$
$$
 \le
2  \int_{Q}| {v}|^2  |{ \nabla\phi ({\frac{x_{3}^{2}}{k^{2}}})}|^2  dx
 +  2  \int_{Q}   |  \phi({\frac{x_{3}^{2}}{k^{2}}})|^2  |{ \nabla v}|^{2} dx 
 $$
 $$
 \leq  2\int_{k \leq |{x_{3}}| \leq \sqrt{2}k} |{v}|^{2}|
{ \phi^\prime ({\frac{x_{3}^{2}}{k^{2}}})}|^{2}
 {\frac{|{2x_{3}}|^{2}}{k^{4}}}dx 
+ 2\int_{Q} \cutf |\nabla v |^2 dx
$$
\be
\label{p43_21}
\leq  {\frac{c}{k^2}} \int_{k \le |x_3| \le \sqrt{2} k} |v|^2 dx 
+ 2 \int_{Q}  \cutf |\nabla v |^2 dx
 \leq  {\frac{c}{k^2}}\| v\|^2+ 2 \int_{Q}  \cutf |\nabla v |^2 dx
\ee
Since $v \in H_{0}^{1}(Q)$ we have $ \phi({\frac{x_{3}^{2}}{k^{2}}}) v \in H_{0}^{1}(Q)$ and hence
by  \eqref{poincare} and \eqref{p43_21} we get
$$
\int_Q \cutf | v |^2 dx \le {\frac 1\lambda} 
\int_Q  | \nabla \left (  \phi({\frac{x_{3}^{2}}{k^{2}}})\  v \right ) |^2  dx
\le {\frac{c}{k^2 \lambda}}\| v\|^2+ {\frac 2\lambda}
 \int_{Q} \cutf |\nabla v |^2 dx,
 $$
 and hence we have
\be
\label{p43_22}
 \int_Q \cutf |\nabla v |^2 dx
 \ge {\frac 12} \lambda   \int_Q \cutf  |v|^2 dx
 -{\frac c{2k^2}} \| v \|^2.
\ee
By \eqref{p43_22} we find that
\be
\label{p43_23}
 2 \nu \int_Q \cutf |\nabla v |^2 dx
 \ge \nu  \int_Q \cutf |\nabla v |^2 dx
 + 
  {\frac 12} \nu \lambda   \int_Q \cutf  |v|^2 dx
 -{\frac {c\nu}{2k^2}} \| v \|^2.
 \ee
 On the other hand, we have
 \be
 \label{p43_24}
 \| v\|_3^3 \le c \| v \|_{H^1}^3
 \le c + \| v \|_{H^1}^4.
 \ee
 By \eqref{p43_20} and \eqref{p43_23}-\eqref{p43_24} we obtain that, for all $k \ge 1$,
 $$
  {\frac d{dt}} \int_Q \cutf (|v|^2 + |\nabla v|^2 ) dx
+   \nu \int_Q \cutf |\nabla v |^2 dx + 
  {\frac 12} \nu \lambda   \int_Q \cutf  |v|^2 dx
$$
$$
\le \left ( {\frac 14} \nu \lambda + \beta |y(\theta_t \omega )| \right )
\int_Q\cutf |v|^2 dx
+ {\frac ck} (\| v \|^2 + \| \nabla v \|^2 ) + {\frac ck} (1 +  \| v\|^4_{H^1}) + {\frac ck} \| \nabla v_t \|^2
$$
$$
+ c\int_Q \cutf \left (
 g^2 + |\zt|^2 + |\De \zt|^2 + |\nabla \zt |^2 + |\zt |^2 \ |\nabla \zt |
\right ) dx
$$
\be
\label{p43_30}
+ {\frac ck} \left ( \| \zt \|^2 + \| \zt \|^3_3 \right ).
\ee
By \eqref{fixbeta}  and \eqref{p43_30} we find that
$$
  {\frac d{dt}} \int_Q \cutf (|v|^2 + |\nabla v|^2 ) dx
+    \int_Q \cutf   \left ( \delta - \beta |y(\theta_t \omega )| \right )  \left ( | v|^2 +  |\nabla v |^2 \right ) dx 
$$
$$
\le   {\frac ck} (\| v \|^2 + \| \nabla v \|^2 ) + {\frac ck} (1 +  \| v\|^4_{H^1}) + {\frac ck} \| \nabla v_t \|^2
$$
$$
+ c\int_Q \cutf \left (
 g^2 + |\zt|^2 + |\De \zt|^2 + |\nabla \zt |^2 + |\zt |^2 \ |\nabla \zt |
\right ) dx
$$
\be
\label{p43_31}
+ {\frac ck} \left ( \| \zt \|^2 + \| \zt \|^3_3 \right ).
\ee
Multiplying  \eqref{p43_31} by $e^{\int_0^t (\delta -\beta |y(\theta_\tau \omega )| ) d \tau}$
and then integrating over $(0,t)$, we get  that, for all $t \ge 0$,
$$
 \int_Q \cutf (|v(t, \omega, v_0(\omega) )|^2 + |\nabla v (t, \omega, v_0(\omega) )  |^2 ) dx
 $$
 $$
 \le  e^{- \int_0^t (\delta -\beta |y(\theta_\tau \omega )| ) d \tau} \int_Q \cutf 
 \left (|v_0(\omega)|^2 + |\nabla v_0(\omega)|^2 \right ) dx
 $$
 $$ +{\frac ck} \int_0^t 
 e^{\int_t^s (\delta -\beta |y(\theta_\tau \omega )| ) d \tau} ds
 +{\frac ck} \int_0^t 
 e^{\int_t^s (\delta -\beta |y(\theta_\tau \omega )| ) d \tau}
 \| v(s, \omega, v_0(\omega) )\|^2_{H^1} ds
 $$
 $$
 +{\frac ck} \int_0^t 
 e^{\int_t^s (\delta -\beta |y(\theta_\tau \omega )| ) d \tau}
 \| v(s, \omega, v_0(\omega) )\|^4_{H^1} ds
 $$
 $$
 +{\frac ck} \int_0^t 
 e^{\int_t^s (\delta -\beta |y(\theta_\tau \omega )| ) d \tau}
 \| \nabla v_s (s, \omega, v_0(\omega) )\|^2ds
 $$
 $$
 +c  \int_0^t 
 e^{\int_t^s (\delta -\beta |y(\theta_\tau \omega )| ) d \tau}
\left (
\int_Q \cutf   \ 
 g^2   dx
\right )ds
 $$
 $$
 +c  \int_0^t 
 e^{\int_t^s (\delta -\beta |y(\theta_\tau \omega )| ) d \tau} 
 \left (
\int_Q \cutf \left (
  |z(\theta_{s} \omega )|^2 + |\De  z(\theta_{s} \omega )|^2 
  \right ) dx \right ) ds
  $$
  $$
 +c  \int_0^t 
 e^{\int_t^s (\delta -\beta |y(\theta_\tau \omega )| ) d \tau} 
 \left (
\int_Q \cutf \left (
  |\nabla  z(\theta_{s} \omega ) |^2 
+ | z(\theta_{s} \omega ) |^2 \ |\nabla  z(\theta_{s} \omega ) |
\right ) dx
\right )
 ds
 $$
 $$
 +{\frac ck} \int_0^t 
 e^{\int_t^s (\delta -\beta |y(\theta_\tau \omega )| ) d \tau}
 \left (
 \| z(\theta_s \omega )\|^2 +  \| z(\theta_s \omega )\|^3_3
 \right ) ds.
 $$
 Replacing $\omega$ by $\theta_{-t} \omega $ in the above, we find that, for all
 $t \ge 0$,
 $$
 \int_Q \cutf (|v(t, \theta_{-t} \omega, v_0(\theta_{-t} \omega) )|^2 
+ |\nabla v (t, \theta_{-t} \omega, v_0(\theta_{-t} \omega) )  |^2 ) dx
 $$
 $$
 \le  e^{- \int_0^t (\delta -\beta |y(\theta_{\tau -t} \omega )| ) d \tau} \int_Q \cutf 
 \left (|v_0(\theta_{-t}\omega)|^2 + |\nabla v_0(\theta_{-t}\omega)|^2 \right ) dx
 $$
 $$ +{\frac ck} \int_0^t 
 e^{\int_t^s (\delta -\beta |y(\theta_{\tau -t} \omega )| ) d \tau} ds
 +{\frac ck} \int_0^t 
 e^{\int_t^s (\delta -\beta |y(\theta_{\tau -t} \omega )| ) d \tau}
 \| v(s, \theta_{-t}\omega, v_0(\theta_{-t}\omega) )\|^2_{H^1} ds
 $$
 $$
 +{\frac ck} \int_0^t 
 e^{\int_t^s (\delta -\beta |y(\theta_{\tau -t} \omega )| ) d \tau}
 \| v(s, \theta_{-t}\omega, v_0(\theta_{-t}\omega) )\|^4_{H^1} ds
 $$
 $$
 +{\frac ck} \int_0^t 
 e^{\int_t^s (\delta -\beta |y(\theta_{\tau-t} \omega )| ) d \tau}
 \| \nabla v_s (s, \theta_{-t}\omega, v_0(\theta_{-t}\omega) )\|^2 ds
 $$
 $$
 +c  \int_0^t 
 e^{\int_t^s (\delta -\beta |y(\theta_{\tau -t} \omega )| ) d \tau}
\left (
\int_Q \cutf   \ 
 g^2   dx
\right )ds
 $$
 $$
 +c  \int_0^t 
 e^{\int_t^s (\delta -\beta |y(\theta_{\tau -t} \omega )| ) d \tau}
 \left (
\int_Q \cutf \left (
  |z(\theta_{s-t} \omega )|^2 + |\De  z(\theta_{s-t} \omega )|^2 
  \right ) dx \right ) ds
  $$
  $$ +c  \int_0^t 
 e^{\int_t^s (\delta -\beta |y(\theta_{\tau -t} \omega )| ) d \tau}
 \left (
\int_Q \cutf \left (
 |\nabla  z(\theta_{s-t} \omega ) |^2 
+ | z(\theta_{s-t} \omega ) |^2 \ |\nabla  z(\theta_{s-t} \omega ) |
\right ) dx \right )
 ds
 $$
 \be
\label{p43_40}
 +{\frac ck} \int_0^t 
 e^{\int_t^s (\delta -\beta |y(\theta_{\tau -t} \omega )| ) d \tau}
 \left (
 \| z(\theta_{s-t} \omega )\|^2 +  \| z(\theta_{s-t} \omega )\|^3_3
 \right ) ds.
\ee
In what follows, we estimate every term on the right-hand side
of \eqref{p43_40}. For the first term,  by \eqref{deltacond} we have
$$
 e^{- \int_0^t (\delta -\beta |y(\theta_{\tau -t} \omega )| ) d \tau} \int_Q \cutf 
 \left (|v_0(\theta_{-t}\omega)|^2 + |\nabla v_0(\theta_{-t}\omega)|^2 \right ) dx
 $$
 $$ \le
  e^{- \delta t  + \beta  \int_0^t    |y(\theta_{\tau -t} \omega )|  d \tau}  
\| v_0(\theta_{-t}\omega)\|^2_{H^1}   
\le
e^{- \delta t  + \beta  \int_{-t}^0    |y(\theta_{\tau } \omega )|  d \tau}  
\| v_0(\theta_{-t}\omega)\|^2_{H^1}   
 $$
\be
\label{p43_41}
 \le
e^{- {\frac 78}\delta t  }
\| v_0(\theta_{-t}\omega)\|^2_{H^1} , \quad \mbox{for all} \ t \ge T_0(\omega).  
 \ee
 Since $v_0(\theta_{-t} \omega)  \in B(\theta_{-t} \omega)$
 and $B=\{B(\omega)\}_{\omega \in \Omega} \in \mathcal{D}$, the right-hand side
 of \eqref{p43_41} tends to zero as $t \to \infty$. 
 Therefore, given $\epsilon>0$, there is $T_1=T_1(B, \omega, \epsilon)>0$
 such that for all $t \ge T_1$,
\be
\label{p43_42}
 e^{- \int_0^t (\delta -\beta |y(\theta_{\tau -t} \omega )| ) d \tau} \int_Q \cutf 
 \left (|v_0(\theta_{-t}\omega)|^2 + |\nabla v_0(\theta_{-t}\omega)|^2 \right ) dx
 \le \epsilon.
 \ee
 Note that the second term 
 on the right-hand side of \eqref{p43_40} satisfies
 $$ 
  {\frac ck} \int_0^t 
 e^{\int_t^s (\delta -\beta |y(\theta_{\tau -t} \omega )| ) d \tau} ds
=
 {\frac ck} \int_0^t 
 e^{\delta (s-t)   -\beta\int_t^s  |y(\theta_{\tau -t} \omega )| ) d \tau} ds
 $$
\be\label{p43_43}
 =
 {\frac ck} \int_0^t 
 e^{\delta (s-t)   -\beta\int_0^{s-t}  |y(\theta_{\tau } \omega )| ) d \tau} ds
 =
 {\frac ck} \int_{-t}^0
 e^{\delta  s   + \beta\int_s^0  |y(\theta_{\tau } \omega )| ) d \tau} ds
\ee
By \eqref{deltacond}, the integrand in \eqref{p43_43} converges to zero
exponentially as $s \to -\infty$, and hence the following integral is well-defined:
\be
\label{p43_44}
r_1(\omega) =
 \int_{-\infty}^0
 e^{\delta  s   + \beta\int_s^0  |y(\theta_{\tau } \omega )| ) d \tau} ds .
\ee
It follows from \eqref{p43_43}-\eqref{p43_44} that, for all $t\ge 0$,
\be
\label{p43_45}
  {\frac ck} \int_0^t 
 e^{\int_t^s (\delta -\beta |y(\theta_{\tau -t} \omega )| ) d \tau} ds
\le 
 {\frac ck} r_1(\omega).
 \ee
 By \eqref{vs}, the third term on the right-hand side of
 \eqref{p43_40} is bounded by
 $$
  {\frac ck} \int_0^t 
 e^{\int_t^s (\delta -\beta |y(\theta_{\tau -t} \omega )| ) d \tau}
 \| v(s, \theta_{-t}\omega, v_0(\theta_{-t}\omega) )\|^2_{H^1} ds
 $$
 $$
\le 
  {\frac ck} \int_0^t 
 e^{\int_t^s (\delta -\beta |y(\theta_{\tau -t} \omega )| ) d \tau}
e^{ -\delta s + \beta \int_{-t}^{s-t}  |y(\theta_{\tau} \omega)|  d \tau} \| v_0 (\theta_{-t}\omega) \|^2_{H^1}
ds
$$
 \be
\label{p43_46}
+  {\frac ck} \int_0^t 
 e^{\int_t^s (\delta -\beta |y(\theta_{\tau -t} \omega )| ) d \tau}
\left (  \int_{-t}^{s-t} (1 + |y(\theta_{\sigma} \omega)|^2 + |y(\theta_{\sigma} \omega)|^3 ) 
e^{  \delta (\sigma -s +t) + \beta \int_\sigma^{s-t}  |y(\theta_{\tau} \omega)|  d \tau} d \sigma
\right ) ds.
\ee
By \eqref{deltacond}, the first term on the right-hand side of \eqref{p43_46}
is given by
$$
 {\frac ck} \int_0^t 
 e^{\int_t^s (\delta -\beta |y(\theta_{\tau -t} \omega )| ) d \tau}
e^{ -\delta s + \beta \int_{-t}^{s-t}  |y(\theta_{\tau} \omega)|  d \tau} \| v_0 (\theta_{-t}\omega) \|^2_{H^1}
ds
$$
$$ = 
  {\frac ck} \int_0^t 
 e^{   - \delta t  +\beta \int_s^t |y(\theta_{\tau -t} \omega )|  d \tau 
 +  \beta \int_{-t}^{s-t}  |y(\theta_{\tau} \omega)|  d \tau } \| v_0 (\theta_{-t}\omega) \|^2_{H^1}
ds
$$
$$ = 
  {\frac ck} \int_0^t 
 e^{   - \delta t  +\beta \int^0_{-t} |y(\theta_{\tau } \omega )|  d \tau }
 \| v_0 (\theta_{-t}\omega) \|^2_{H^1}
ds
$$
 $$ = {\frac ck} 
    t 
 e^{   - \delta t  +\beta \int^0_{-t} |y(\theta_{\tau } \omega )|  d \tau }
 \| v_0 (\theta_{-t}\omega) \|^2_{H^1} 
\le 
{\frac ck}  
    t 
 e^{   -{\frac 78}  \delta t }  
 \| v_0 (\theta_{-t}\omega) \|^2_{H^1} ,
$$
for all $t \ge T_0(\omega)$. Since  $t 
 e^{   -{\frac 78}  \delta t }  
 \| v_0 (\theta_{-t}\omega) \|^2_{H^1}$
 tends to zero  as $t \to \infty$, there is $T_2=T_2(B,\omega)>0$ such that
 for all $t \ge T_2$,
\be
\label{p43_50}
 {\frac ck} \int_0^t 
 e^{\int_t^s (\delta -\beta |y(\theta_{\tau -t} \omega )| ) d \tau}
e^{ -\delta s + \beta \int_{-t}^{s-t}  |y(\theta_{\tau} \omega)|  d \tau} \| v_0 (\theta_{-t}\omega) \|^2_{H^1}
ds
\le {\frac ck}.
\ee
For the second term on the right-hand side of
\eqref{p43_46} we have
$$
 {\frac ck} \int_0^t 
 e^{\int_t^s (\delta -\beta |y(\theta_{\tau -t} \omega )| ) d \tau}
\left (  \int_{-t}^{s-t} (1 + |y(\theta_{\sigma} \omega)|^2 + |y(\theta_{\sigma} \omega)|^3 ) 
e^{  \delta (\sigma -s +t) + \beta \int_\sigma^{s-t}  |y(\theta_{\tau} \omega)|  d \tau} d \sigma
\right ) ds
 $$
 $$ \le
 {\frac ck} \int_0^t 
 e^{\delta (s-t)  + \beta \int_s^t |y(\theta_{\tau -t} \omega )|  d \tau}
\left (  \int_{-t}^{s-t} (1 + |y(\theta_{\sigma} \omega)|^2 + |y(\theta_{\sigma} \omega)|^3 ) 
e^{  {\frac 12}\delta (\sigma -s +t) + \beta \int_\sigma^{s-t}  |y(\theta_{\tau} \omega)|  d \tau} d \sigma
\right ) ds
 $$
 $$ \le
 {\frac ck} \int_0^t 
 e^{{\frac 12} \delta (s-t)  + \beta \int_s^t |y(\theta_{\tau -t} \omega )|  d \tau}
\left (  \int_{-t}^{s-t} (1 + |y(\theta_{\sigma} \omega)|^2 + |y(\theta_{\sigma} \omega)|^3 ) 
e^{  {\frac 12}\delta \sigma  + \beta \int_\sigma^{s-t}  |y(\theta_{\tau} \omega)|  d \tau} d \sigma
\right ) ds
 $$
 $$ \le
 {\frac ck} \int_0^t 
 e^{{\frac 12} \delta (s-t)  + \beta \int_s^t |y(\theta_{\tau -t} \omega )|  d \tau} ds
 \int_{-t}^{0} (1 + |y(\theta_{\sigma} \omega)|^2 + |y(\theta_{\sigma} \omega)|^3 ) 
e^{  {\frac 12}\delta \sigma  + \beta \int_\sigma^{0}  |y(\theta_{\tau} \omega)|  d \tau} d \sigma
 $$ 
  \be
  \label{p43_51} \le
 {\frac ck} 
 \int^0_{-t}
 e^{{\frac 12} \delta s  + \beta \int_s^0 |y(\theta_{\tau } \omega )|  d \tau} ds
 \int_{-t}^{0} (1 + |y(\theta_{\sigma} \omega)|^2 + |y(\theta_{\sigma} \omega)|^3 ) 
e^{  {\frac 12}\delta \sigma  + \beta \int_\sigma^{0}  |y(\theta_{\tau} \omega)|  d \tau} d \sigma
\ee
By \eqref{deltacond} we know that the following integrals are convergent:
\be
\label{p43_52}
r_2 (\omega)
= \int^0_{-\infty}
 e^{{\frac 12} \delta s  + \beta \int_s^0 |y(\theta_{\tau } \omega )|  d \tau} ds,
\ee 
and
\be
\label{p43_53} r_3(\omega) = 
\int_{-\infty}^{0} (1 + |y(\theta_{\sigma} \omega)|^2 + |y(\theta_{\sigma} \omega)|^3 ) 
e^{  {\frac 12}\delta \sigma  + \beta \int_\sigma^{0}  |y(\theta_{\tau} \omega)|  d \tau} d \sigma.
\ee 
By \eqref{p43_51}-\eqref{p43_53} we obtain that
$$
 {\frac ck} \int_0^t 
 e^{\int_t^s (\delta -\beta |y(\theta_{\tau -t} \omega )| ) d \tau}
\left (  \int_{-t}^{s-t} (1 + |y(\theta_{\sigma} \omega)|^2 + |y(\theta_{\sigma} \omega)|^3 ) 
e^{  \delta (\sigma -s +t) + \beta \int_\sigma^{s-t}  |y(\theta_{\tau} \omega)|  d \tau} d \sigma
\right ) ds
$$
\be
\label{p43_60}
\le {\frac ck} r_2(\omega) r_3(\omega).
\ee
Then it follows  from \eqref{p43_46} and \eqref{p43_60} that, for
all $t \ge T_2$,
\be
\label{p43_61}
  {\frac ck} \int_0^t 
 e^{\int_t^s (\delta -\beta |y(\theta_{\tau -t} \omega )| ) d \tau}
 \| v(s, \theta_{-t}\omega, v_0(\theta_{-t}\omega) )\|^2_{H^1} ds
 \le {\frac ck} (1 + r_1(\omega) r_2(\omega) ).
\ee

For the fourth term on the right-hand side of \eqref{p43_40},
by \eqref{vs},  we have
 $$
  {\frac ck} \int_0^t 
 e^{\int_t^s (\delta -\beta |y(\theta_{\tau -t} \omega )| ) d \tau}
 \| v(s, \theta_{-t}\omega, v_0(\theta_{-t}\omega) )\|^4_{H^1} ds
 $$
 $$
\le 
  {\frac ck} \int_0^t 
 e^{\int_t^s (\delta -\beta |y(\theta_{\tau -t} \omega )| ) d \tau}
\left (
e^{ -\delta s + \beta \int_{-t}^{s-t}  |y(\theta_{\tau} \omega)|  d \tau} \| v_0 (\theta_{-t}\omega) \|^2_{H^1}
\right )^2
ds
$$
 \be
\label{p43_62}
+  {\frac ck} \int_0^t 
 e^{\int_t^s (\delta -\beta |y(\theta_{\tau -t} \omega )| ) d \tau}
\left (  \int_{-t}^{s-t} (1 + |y(\theta_{\sigma} \omega)|^2 + |y(\theta_{\sigma} \omega)|^3 ) 
e^{  \delta (\sigma -s +t) + \beta \int_\sigma^{s-t}  |y(\theta_{\tau} \omega)|  d \tau} d \sigma
\right )^2 ds.
\ee
We now deal with the first term on the right-hand side of
the above, which is given by
$$
  {\frac ck} \int_0^t 
 e^{\int_t^s (\delta -\beta |y(\theta_{\tau -t} \omega )| ) d \tau}
\left (
e^{ -\delta s + \beta \int_{-t}^{s-t}  |y(\theta_{\tau} \omega)|  d \tau} \| v_0 (\theta_{-t}\omega) \|^2_{H^1}
\right )^2
ds
$$
$$
=
  {\frac ck} \int_0^t 
 e^{ -\delta t - \delta s  + 2 \beta \int_{-t}^{s-t}
  |y(\theta_{\tau } \omega )|  d \tau
+  \beta \int_{s-t}^{0}  |y(\theta_{\tau} \omega)|  d \tau}
 \| v_0 (\theta_{-t}\omega) \|^4_{H^1}
ds
$$
   $$
=
  {\frac ck} \int_0^t 
 e^{ -\delta t - \delta s  +  \beta \int_{-t}^{0}
  |y(\theta_{\tau } \omega )|  d \tau
+  \beta \int_{-t}^{s-t}  |y(\theta_{\tau} \omega)|  d \tau}
 \| v_0 (\theta_{-t}\omega) \|^4_{H^1}
ds
$$    
$$
\le
  {\frac ck} \int_0^t 
 e^{ -\delta t - \delta s  + 2 \beta \int_{-t}^{0}
  |y(\theta_{\tau } \omega )|  d \tau}
 \| v_0 (\theta_{-t}\omega) \|^4_{H^1}
ds   
$$
\be
\label{p43_63}
\le
  {\frac c{k\delta}} 
 e^{ -\delta t  + 2 \beta \int_{-t}^{0}
  |y(\theta_{\tau } \omega )|  d \tau}
 \| v_0 (\theta_{-t}\omega) \|^4_{H^1}.
\ee
By \eqref{deltacond}, we know that
 $$
 e^{ -\delta t  + 2 \beta \int_{-t}^{0}
  |y(\theta_{\tau } \omega )|  d \tau}
 \| v_0 (\theta_{-t}\omega) \|^4_{H^1}  
 \to 0 \quad \mbox{as} \quad t \to \infty,
 $$ and hence it follows from \eqref{p43_63}
 that there is $T_3=T_3(B, \omega)$ such that
 for all $t \ge T_3$,
\be
\label{p43_64}
  {\frac ck} \int_0^t 
 e^{\int_t^s (\delta -\beta |y(\theta_{\tau -t} \omega )| ) d \tau}
\left (
e^{ -\delta s + \beta \int_{-t}^{s-t}  |y(\theta_{\tau} \omega)|  d \tau} \| v_0 (\theta_{-t}\omega) \|^2_{H^1}
\right )^2
ds
\le {\frac ck}.
\ee
For the second term on the right-hand side of \eqref{p43_62}, we
have
$$
  {\frac ck} \int_0^t 
 e^{\int_t^s (\delta -\beta |y(\theta_{\tau -t} \omega )| ) d \tau}
\left (  \int_{-t}^{s-t} (1 + |y(\theta_{\sigma} \omega)|^2 + |y(\theta_{\sigma} \omega)|^3 ) 
e^{  \delta (\sigma -s +t) + \beta \int_\sigma^{s-t}  |y(\theta_{\tau} \omega)|  d \tau} d \sigma
\right )^2 ds
$$
$$
\le   {\frac ck} \int_0^t 
 e^{  \delta (s-t)   +\beta \int_s^t  |y(\theta_{\tau -t} \omega )|  d \tau}
\left (  \int_{-t}^{s-t} (1 + |y(\theta_{\sigma} \omega)|^2 + |y(\theta_{\sigma} \omega)|^3 ) 
e^{  {\frac 38} \delta (\sigma -s +t) + \beta \int_\sigma^{s-t}  |y(\theta_{\tau} \omega)|  d \tau} d \sigma
\right )^2 ds
$$
 $$
\le   {\frac ck} \int_0^t 
 e^{ {\frac 14} \delta (s-t)   + \beta \int_s^t |y(\theta_{\tau -t} \omega )|  d \tau}
\left (  \int_{-t}^{s-t} (1 + |y(\theta_{\sigma} \omega)|^2 + |y(\theta_{\sigma} \omega)|^3 ) 
e^{  {\frac 38} \delta \sigma + \beta \int_\sigma^{s-t}  |y(\theta_{\tau} \omega)|  d \tau} d \sigma
\right )^2 ds
$$
$$
\le   {\frac ck} \int_0^t 
 e^{ {\frac 14} \delta (s-t)   + \beta\int^0_{s-t}  |y(\theta_{\tau } \omega )|  d \tau} ds
\left (  \int_{-t}^{0} (1 + |y(\theta_{\sigma} \omega)|^2 + |y(\theta_{\sigma} \omega)|^3 ) 
e^{  {\frac 38} \delta \sigma + \beta \int_\sigma^{0}  |y(\theta_{\tau} \omega)|  
d \tau} d \sigma
\right )^2 
$$
\be
\label{p43_65}
\le   {\frac ck} \int_{-t}^0
 e^{ {\frac 14} \delta s    + \beta\int^0_{s}  |y(\theta_{\tau } \omega )|  d \tau} ds
\left (  \int_{-t}^{0} (1 + |y(\theta_{\sigma} \omega)|^2 + |y(\theta_{\sigma} \omega)|^3 ) 
e^{  {\frac 38} \delta \sigma + \beta \int_\sigma^{0}  |y(\theta_{\tau} \omega)|  
d \tau} d \sigma
\right )^2 .
\ee
Note that \eqref{deltacond} implies the convergence of the integrals:
\be
\label{p43_66}
r_4(\omega) = 
\int_{-\infty}^0
 e^{ {\frac 14} \delta s    + \beta\int^0_{s}  |y(\theta_{\tau } \omega )|  d \tau} ds,
 \ee
 and
 \be
\label{p43_67}
r_5(\omega) = 
 \int_{-\infty}^{0} (1 + |y(\theta_{\sigma} \omega)|^2 + |y(\theta_{\sigma} \omega)|^3 ) 
e^{  {\frac 38} \delta \sigma + \beta \int_\sigma^{0}  |y(\theta_{\tau} \omega)|  
d \tau} d \sigma.
 \ee
 Therefore it follows from \eqref{p43_65}-\eqref{p43_67} that, for all $t\ge 0$,
 $$
  {\frac ck} \int_0^t 
 e^{\int_t^s (\delta -\beta |y(\theta_{\tau -t} \omega )| ) d \tau}
\left (  \int_{-t}^{s-t} (1 + |y(\theta_{\sigma} \omega)|^2 + |y(\theta_{\sigma} \omega)|^3 ) 
e^{  \delta (\sigma -s +t) + \beta \int_\sigma^{s-t}  |y(\theta_{\tau} \omega)|  d \tau} d \sigma
\right )^2 ds
$$
\be
\label{p43_70}
\le
{\frac ck} r_4(\omega) \  r_5^2(\omega).
\ee
By \eqref{p43_62}, \eqref{p43_64} and \eqref{p43_70} we find that,
for all $t \ge T_3$,
\be\label{p43_71}
  {\frac ck} \int_0^t 
 e^{\int_t^s (\delta -\beta |y(\theta_{\tau -t} \omega )| ) d \tau}
 \| v(s, \theta_{-t}\omega, v_0(\theta_{-t}\omega) )\|^4_{H^1} ds
 \le {\frac ck} (1 + r_4(\omega) \ r_5^2(\omega) ).
\ee
Note that $g \in \h$, and hence for given $\epsilon>0$, there is
$k_1 =k_1(\epsilon)>0$ such that, for all $k\ge k_1$,
$$
\int_{|x_3| \ge k} g^2(x) dx \le \epsilon,
$$
from which the sixth term on the right-hand side of \eqref{p43_40} is
bounded by
$$
  c  \int_0^t 
 e^{\int_t^s (\delta -\beta |y(\theta_{\tau -t} \omega )| ) d \tau}
\left (
\int_Q \cutf   \ 
 g^2   dx
\right )ds
 $$
 $$ \le 
  c  \int_0^t 
 e^{\int_t^s (\delta -\beta |y(\theta_{\tau -t} \omega )| ) d \tau}
\left (
\int_{|x_3| \ge k}  \cutf   \ 
 g^2   dx
\right )ds
 $$
\be
\label{p43_72}\le 
  \epsilon c  \int_0^t 
 e^{\int_t^s (\delta -\beta |y(\theta_{\tau -t} \omega )| ) d \tau}
 ds
 \le \epsilon c \  r_1(\omega),
 \ee
 where $r_1(\omega)$ is given by
 \eqref{p43_44}.
 Since $(I-\De)^{-1} h \in H^2(Q)$, there is $k_2 =k_2(\omega)>0$ such that
 for all $k \ge k_2$,
 \be
 \label{p43_73}
 \int_{|x_3| \ge k}    \left (
|(I-\De )^{-1} h |^2  +| \De (I-\De )^{-1} h |^2
\right ) dx
\le \epsilon.
\ee
Note that $\zt = (I-\De)^{-1} h y(\theta_t \omega)$. By 
\eqref{p43_73} the seventh term on the right-hand side
of \eqref{p43_40} satisfies
$$
  c  \int_0^t 
 e^{\int_t^s (\delta -\beta |y(\theta_{\tau -t} \omega )| ) d \tau}
 \left (
\int_Q \cutf \left (
  |z(\theta_{s-t} \omega )|^2 + |\De  z(\theta_{s-t} \omega )|^2 
  \right ) dx \right ) ds
  $$
  $$
  \le c  \int_0^t 
 e^{\int_t^s (\delta -\beta |y(\theta_{\tau -t} \omega )| ) d \tau}
 \left (
\int_{|x_3| \ge k} \cutf \left (
  |z(\theta_{s-t} \omega )|^2 + |\De  z(\theta_{s-t} \omega )|^2 
  \right ) dx \right ) ds
  $$
  $$
  \le  \epsilon c  \int_0^t 
 e^{\int_t^s (\delta -\beta |y(\theta_{\tau -t} \omega )| ) d \tau}
  |y(\theta_{s-t}\omega)|^2 ds
  \le  \epsilon c  \int_0^t 
 e^{  \delta (s-t)  + \beta \int_{s-t}^0  |y(\theta_{\tau } \omega )|  d \tau}
  |y(\theta_{s-t}\omega)|^2 ds
  $$
\be
\label{p43_74}
   \le  \epsilon c  \int_{-t}^0
 e^{  \delta  s + \beta \int_{s}^0  |y(\theta_{\tau } \omega )|  d \tau}
  |y(\theta_{s}\omega)|^2 ds
  \le \epsilon c \  r_6(\omega),
 \ee
  where   $r_6(\omega)$ is given by
  $$ r_6(\omega) = 
  \int_{-\infty}^0
 e^{  \delta  s + \beta \int_{s}^0  |y(\theta_{\tau } \omega )|  d \tau}
  |y(\theta_{s}\omega)|^2 ds.
  $$
  Note that $r_6(\omega)$ is well-defined  by \eqref{deltacond}.
  Similarly,  we  can find a random function  $r_7(\omega)$
 such that the eighth term on the right-hand
  side of \eqref{p43_40}  satisfies 
  \be
\label{p43_75} 
 c  \int_0^t 
 e^{\int_t^s (\delta -\beta |y(\theta_{\tau -t} \omega )| ) d \tau}
 \left (
\int_Q \cutf \left (
 |\nabla  z(\theta_{s-t} \omega ) |^2 
+ | z(\theta_{s-t} \omega ) |^2 \ |\nabla  z(\theta_{s-t} \omega ) |
\right ) dx \right )
 ds  \le \epsilon c  \  r_7(\omega).
\ee
For the last term on the right-hand side of \eqref{p43_40} we have
$$
  {\frac ck} \int_0^t 
 e^{\int_t^s (\delta -\beta |y(\theta_{\tau -t} \omega )| ) d \tau}
 \left (
 \| z(\theta_{s-t} \omega )\|^2 +  \| z(\theta_{s-t} \omega )\|^3_3
 \right ) ds
$$
$$ \le
  {\frac ck} \int_0^t 
 e^{\int_t^s (\delta -\beta |y(\theta_{\tau -t} \omega )| ) d \tau}
 \left (
 | y(\theta_{s-t} \omega )|^2 +  | y(\theta_{s-t} \omega )|^3
 \right ) ds
$$
$$ 
  \le {\frac ck}   \int_0^t 
 e^{  \delta (s-t)  + \beta \int_{s-t}^0   |y(\theta_{\tau } \omega )|  d \tau}
 (  |y(\theta_{s-t}\omega)|^2  + |y(\theta_{s-t}\omega)|^3 )ds
  $$
\be
\label{p43_74}
   \le {\frac ck}   \int_{-t}^0
 e^{  \delta  s + \beta \int_{s}^0    |y( \theta_{\tau } \omega )| 
 ) d \tau}
(   |y(\theta_{s}\omega)|^2  + |y(\theta_{s}\omega)|^3 )ds
  \le {\frac ck}  \  r_8(\omega),
 \ee
 where $r_8(\omega)$ is given by
 $$
 r_8(\omega)
 =  \int_{-\infty}^0
 e^{  \delta  s + \beta \int_{s}^0    |y( \theta_{\tau } \omega )| 
 ) d \tau}
(   |y(\theta_{s}\omega)|^2  + |y(\theta_{s}\omega)|^3 )ds
 .$$
 We now deal with the fifth term on the right-hand side 
 of \eqref{p43_40}. By Lemma \ref{lem42} we have
$$
  {\frac ck} \int_0^t 
 e^{\int_t^s (\delta -\beta |y(\theta_{\tau-t} \omega )| ) d \tau}
 \| \nabla v_s (s, \theta_{-t}\omega, v_0(\theta_{-t}\omega) )\|^2 ds
 $$
$$
\le {\frac ck} \int_0^t 
 e^{\int_t^s (\delta -\beta |y(\theta_{\tau-t} \omega )| ) d \tau} ds
 $$
 $$
  +
{\frac ck} \int_0^t 
 e^{\int_t^s (\delta -\beta |y(\theta_{\tau-t} \omega )| ) d \tau}
   e^{ - 2 \delta s + 2 \beta \int_{-t}^{s-t}  |y(\theta_{\tau} \omega)|  d \tau}
 \| v_0 (\theta_{-t}\omega) \|^4_{H^1} ds
$$
$$+ 
 {\frac ck} \int_0^t 
 e^{\int_t^s (\delta -\beta |y(\theta_{\tau-t} \omega )| ) d \tau}
 \left ( \int_{-t}^{s-t} (1 + |y(\theta_{\sigma} \omega)|^2 + |y(\theta_{\sigma} \omega)|^3 ) 
e^{  \delta (\sigma -s +t) + \beta \int_\sigma^{s-t}  |y(\theta_{\tau} \omega)|  d \tau} d \sigma
\right )^2 ds
$$
\be
\label{p43_80}
+
 {\frac ck} \int_0^t 
 e^{\int_t^s (\delta -\beta |y(\theta_{\tau-t} \omega )| ) d \tau}
 \left (
\| z(\theta_{s-t} \omega)\|^2
+ \| z(\theta_{s-t} \omega)\|^4_{H^1}
+\| z(\theta_{s-t} \omega)\|^2_{H^2}
\right ) ds.
\ee
Note that the last term of the above has estimates similar to 
\eqref{p43_74}, which along with \eqref{p43_45}, \eqref{p43_64} and
\eqref{p43_70} imply that, there are   $r_9(\omega)$
and  $T_4 =T_4(B,\omega)>0$ such that for all $t \ge T_4$,
\be
\label{p43_81}
  {\frac ck} \int_0^t 
 e^{\int_t^s (\delta -\beta |y(\theta_{\tau-t} \omega )| ) d \tau}
 \| \nabla v_s (s, \theta_{-t}\omega, v_0(\theta_{-t}\omega) )\|^2 ds
 \le {\frac ck} (1 + r_9(\omega ) ).
\ee
Let $T_5=T_5(B,\omega, \epsilon) =\max\{T_1, T_2, T_3, T_4\}$
and $
k_3 =k_3(\epsilon)=\max\{k_1, k_2\}$.
Then it follows from \eqref{p43_40}, \eqref{p43_42},
 \eqref{p43_45}, \eqref{p43_61}, \eqref{p43_71}-\eqref{p43_74}
 and \eqref{p43_81} that, for all $t\ge T_5$ and $k\ge k_3$,
$$ 
 \int_Q \cutf (|v(t, \theta_{-t} \omega, v_0(\theta_{-t} \omega) )|^2 
+ |\nabla v (t, \theta_{-t} \omega, v_0(\theta_{-t} \omega) )  |^2 ) dx
$$
\be\label{p43_90}
\le \epsilon (1 + r_{10} (\omega) ) + {\frac ck} r_{10} (\omega),
\ee
 where $r_{10} (\omega)$ is a positive random function.
By  \eqref{p43_90} we find that there is $k_4=k_4(\omega, \epsilon)>0$
such that for all $t\ge T_5$ and $k\ge k_4$,
$$
 \int_{|x_3|\ge \sqrt{2} k} (|v(t, \theta_{-t} \omega, v_0(\theta_{-t} \omega) )|^2 
+ |\nabla v (t, \theta_{-t} \omega, v_0(\theta_{-t} \omega) )  |^2 ) dx
$$
$$ \le 
 \int_Q \cutf (|v(t, \theta_{-t} \omega, v_0(\theta_{-t} \omega) )|^2 
+ |\nabla v (t, \theta_{-t} \omega, v_0(\theta_{-t} \omega) )  |^2 ) dx
\le \epsilon (2 + r_{10} (\omega) ) ,
$$
which completes the proof.
\end{proof}

In the sequel, we derive uniform estimates
of the solutions on bounded domains which are necessary
for verifying the asymptotic compactness of the stochastic \bbm
equation.
   To this end, we define
  $ \psi =1-\phi $ where $\phi$ is  the function given in \eqref{phif}. 
 Fix $ k \ge 1 $ and
let $ \vb (x, t, \omega )=\psi (\frac{x_{3}^{2}}{k^{2}}) v(x, t, \omega) $.
Then  
 $ \vb(\cdot, t, \omega)  \in H_{0}^{1}(Q_{2k}) $ and
\be
\label{vbv}
 \| \vb(t, \omega)\|_{H^1_0(Q_{2k})}
 \le c 
 \|  v (t, \omega)\|_{H^1(Q_{2k})},
 \quad \forall t \ge 0, \ \omega \in \Omega,
\ee
 where $c$ is a positive deterministic constant,  independent of 
 $\omega \in \Omega$ and $k\ge 1$.
 Note   that
\be
\label{vb1}
\vb_{t}=\psi v_{t},
\ee
\begin{equation}
\label{vb2}
\Delta \vb =(\Delta \psi )v +2\nabla \psi \cdot \nabla v +\psi \Delta v,
\end{equation}
\begin{equation}
\label{vb3}
\Delta \vb_{t}=(\Delta \psi ) v_{t}+2\nabla \psi \cdot \nabla v_{t}+\psi
\Delta v_{t}.
\end{equation}
By \eqref{vb2} and \eqref{vb3}  we have
\begin{equation}
\label{vb4}
\psi \Delta v=\Delta \vb -v\Delta \psi -2\nabla \psi \cdot  \nabla v,
\end{equation}
and
\be
\label{vb5}
 \psi \Delta v_{t}=\Delta \vb_{t}-v_{t}\Delta \psi -2\nabla \psi  \cdot \nabla
v_{t} .
\ee
Multiplying  \eqref{ve}  by $ \psi $ we  get that
\be
\label{vb6}
\psi v_t -  \psi \De v_t - \nu  \psi \De v =
 - \psi \nabla \cdot \vecf (v+z(\theta_t \omega)) 
+ \psi g + \alpha  \psi z(\theta_t \omega) + (\nu -\alpha)
\psi  \De z(\theta_t \omega).
\ee 
Substituting \eqref{vb1} and \eqref{vb4}-\eqref{vb5}  into \eqref{vb6} 
 we find that
$$
\vb _{t}-\Delta \vb_{t}-\nu \Delta  \vb =  
- \psi \nabla \cdot \vecf (v+z(\theta_t \omega)) 
+ \psi g + \alpha  \psi z(\theta_t \omega) + (\nu -\alpha)
\psi  \De z(\theta_t \omega) 
$$
\be 
\label{vbe}
- v_t \De \psi -2 \nabla \psi \cdot \nabla v_t
-\nu v \De \psi -2 \nu \nabla \psi \cdot \nabla v.
\ee
Consider the eigenvalue problem:
\be 
\label{eigen}
-\Delta \vb=\lambda \vb \quad  {\mbox{in}}  \quad
 Q_{2k}, \quad  \vb |_{\partial Q_{2k}}=0.
\ee
Then problem
\eqref{eigen}  has a family of eigenfunctions $\{e_{j}\}_{j=1}^{\infty}$
with corresponding  eigenvalues $\{\lambda_j\}_{j=1}^\infty$
such that
$\{e
_{j}\}_{j=1}^{\infty}$ is an
  orthonormal basis of $L^{2}(Q_{2k})$ and
$$ \lambda _{1}\leq \lambda _{2}\leq \cdots \leq \lambda
_{j}\rightarrow \infty \quad {\mbox{as}}  \  j\rightarrow \infty . $$
Given  $n$, let $X_{n} =$ span$\{e_{1}, \cdots , e_{n}\}$ and $P_{n}:
L^{2}(Q_{2k})\rightarrow X_{n}$ be the projection operator.
For $\vb$, we have the following estimates in $H^1_0(Q_{2k})$.

\begin{lem}
\label{lem44}
 Assume that $g \in \h$, $ h  \in  H^1_0 (Q)$ and
 \eqref{fcond} holds. Let
 $B=\{B(\omega)\}_{\omega \in \Omega}\in \mathcal{D}$ and
 $v_0(\omega) \in B(\omega)$. Then for 
every $\epsilon>0$   and $P$-a.e. $\omega \in \Omega$,
  there exist $T=T (B, \omega, \epsilon)>0$ and $N=N(\omega, \epsilon)>0$
  such that for all  $k\ge 1$,  $t \ge T$ and $n \ge N$,
  $$
  \| (I-P_n) \vb (t, \theta_{-t} \omega, \vb_0(\theta_{-t} \omega ) )\|_{H^1_0(Q_{2k})}
  \le \epsilon.
  $$
 \end{lem}
 
 \begin{proof}
 Let  $ \vb_{n,1}=P_{n}\vb $ and $ \vb_{n,2}= \vb -\vb_{n,1}
$. 
Then applying $I-P_{n}$ to \eqref{vbe}  and   taking the inner product of the 
resulting equation with $\vb_{n,2}$ in $L^{2}(Q_{2k})$  we obtain that
$$
{\frac 12}
 {\frac d{dt}}  
\left ( \| \vb_{n,2} \|^2
+ \| \nabla \vb_{n,2} \|^2
\right ) + \nu \| \nabla \vb_{n,2} \|^2
= 
-  \left ( \psi \nabla \cdot \vecf (v+z(\theta_t \omega)) , \vb_{n,2} \right )
+ (\psi g,  \vb_{n,2})
$$
\be
\label{p44_1}
 + \left ( \alpha  \psi z(\theta_t \omega) + (\nu -\alpha)
\psi  \De z(\theta_t \omega), \vb_{n,2} \right ) 
-
 \left ( v_t \De \psi +2 \nabla \psi \cdot \nabla v_t
+\nu v \De \psi +2 \nu \nabla \psi \cdot \nabla v, \vb_{n,2} \right ).
\ee
By \eqref{fcond}, the nonlinear term in the above is bounded 
by
$$
| \left ( \psi \nabla \cdot \vecf (v+z(\theta_t \omega)) , \vb_{n,2} \right )|
\le \int_{Q_{2k}} \psi \left ( {\frac {x_3^2}{k^2}} \right )
|\vecf^\prime (v+z(\theta_t \omega)) | \ |\nabla v + \nabla z(\theta_t \omega)| 
\ |\vb_{n,2}| dx
$$
$$
\le \int_{Q_{2k}} 
| \gamma_1 + \gamma_2  (v+z(\theta_t \omega)) | \ |\nabla v + \nabla z(\theta_t \omega)| 
\ |\vb_{n,2}| dx
$$
$$
\le
c \left ( \| \nabla v \| + \|\nabla \zt \| \right )
\| \vb_{n,2} \|
+ c  \left ( \| \nabla v \| + \|\nabla \zt \| \right )
\left (\| v\|_6 + \| \zt \|_6 \right )
\| \vb_{n,2} \|_3
$$
$$
\le
c \left ( \|  v \|_{H^1}  + \|  \zt \|_{H^1} \right )
\| \vb_{n,2} \|
+ c  \left ( \|  v \|_{H^1} + \|  \zt \| _{H^1} \right )
\left (\| v\|_{H^1} + \| \zt \|_{H^1} \right )
\| \vb_{n,2} \|_3
$$
$$
\le
c \left ( \|  v \|_{H^1}  + \|  \zt \|_{H^1} \right )
\| \vb_{n,2} \|
+ c  \left ( \|  v \|_{H^1}^2 + \|  \zt \| _{H^1}^2 \right )
\| \nabla \vb_{n,2} \|^{\frac 12} \| \vb_{n,2} \|^{\frac 12}
$$
$$
\le
c \lambda_{n+1}^{-\frac 12} \left ( \|  v \|_{H^1}  + \|  \zt \|_{H^1} \right )
\| \nabla  \vb_{n,2} \|
+ c \lambda_{n+1}^{-\frac 14} 
  \left ( \|  v \|_{H^1}^2 + \|  \zt \| _{H^1}^2 \right )
\| \nabla \vb_{n,2} \| 
$$
$$
\le {\frac 1{16}} \nu \| \nabla \vb_{n,2} \|^2
+ 
c \lambda_{n+1}^{-1} \left ( \|  v \|_{H^1} ^2 + \|  \zt \|_{H^1}^2 \right )
+ c \lambda_{n+1}^{-\frac 12} 
  \left ( \|  v \|_{H^1}^4 + \|  \zt \| _{H^1}^4 \right )
$$
$$
\le {\frac 1{16}} \nu \| \nabla \vb_{n,2} \|^2
+ 
c \lambda_{n+1}^{-1}  
+ c \left (  \lambda_{n+1}^{-\frac 12} + \lambda_{n+1}^{-1} \right )
  \left ( \|  v \|_{H^1}^4 + \|  \zt \| _{H^1}^4 \right ).
$$
\be
\label{p44_2}
\le {\frac 1{16}} \nu \| \nabla \vb_{n,2} \|^2
+ 
c \lambda_{n+1}^{-1}  
+ c \left (  \lambda_{n+1}^{-\frac 12} + \lambda_{n+1}^{-1} \right )
  \left ( \|  v \|_{H^1}^4 + |y(\theta_t \omega) |^4 \right ).
\ee
Note that the second term  on the right-hand side of 
\eqref{p44_1} is bounded by
\be
\label{p44_3}
|(\psi g, \vb_{n,2} ) |
\le \| g \| \ \|  \vb_{n,2} \|
\le \lambda_{n+1}^{-\frac 12} \| g \| \  \| \nabla \vb_{n,2} \|
\le {\frac 1{16} } \nu \| \nabla \vb_{n,2} \|^2
+ c \lambda_{n+1}^{-1} .
\ee
For the third term on the right-hand side of \eqref{p44_1} we have
$$
|\left ( \alpha  \psi z(\theta_t \omega) + (\nu -\alpha)
\psi  \De z(\theta_t \omega), \vb_{n,2} \right ) |
\le c ( \| \zt \| + \De \zt \| )   \|  \vb_{n,2} \|
$$
\be
\label{p44_4}
\le c |y(\theta_t \omega)| \|  \vb_{n,2} \|
\le c \lambda_{n+1}^{-\frac 12} |y(\theta_t \omega)| \ \| \nabla \vb_{n,2} \|
\le
{\frac 1{16}} \nu \| \nabla \vb_{n,2} \|^2
+ c \lambda_{n+1}^{-1}  |y(\theta_t \omega)|^2.
\ee
Similarly, we can check that the last term on the right-hand side 
of \eqref{p44_1} is bounded by
$$
|\left ( v_t \De \psi +2 \nabla \psi \cdot \nabla v_t
+\nu v \De \psi +2 \nu \nabla \psi \cdot \nabla v, \vb_{n,2} \right )|
$$
$$
\le 
{\frac 1{16}} \nu \| \nabla \vb_{n,2} \|^2
+ c \lambda_{n+1}^{-1}  \left (
\| v \|_{H^1}^2 + \| v_t \|_{H^1}^2 \right )
$$
\be
\label{p44_5}
\le 
{\frac 1{16}} \nu \| \nabla \vb_{n,2} \|^2
+ c \lambda_{n+1}^{-1}  \left (1+
\| v \|_{H^1}^4 + \| v_t \|_{H^1}^2 \right ).
\ee
Then it follows from \eqref{p44_1}-\eqref{p44_5} that
$$
{\frac d{dt}} \| \vb_{n,2}\|^2_{H^1}
+ {\frac 32} \nu \| \nabla \vb_{n,2} \|^2
\le 
c \left (  \lambda_{n+1}^{-\frac 12} + \lambda_{n+1}^{-1} \right )
  \left ( \|  v \|_{H^1}^4 + |y(\theta_t \omega) |^4 \right )
  $$
  \be
  \label{p44_6}
 +
  c \lambda_{n+1}^{-1}   \left (1+
\| v \|_{H^1}^4 + \| v_t \|_{H^1}^2 +  |y(\theta_t \omega) |^2  \right ).
\ee
Given  $\epsilon>0$, take $N=N(\epsilon)>0$ large enough such that
for all $n \ge N$, 
\be
\label{p44_10}
\lambda_{n+1} \ge \max\{1, \lambda\}  \quad \mbox{and}
\quad \lambda_{n+1}^{-\frac 12} \le \epsilon,
\ee
where $\lambda$ is the positive constant in\eqref{poincare}.
By \eqref{p44_6} and \eqref{p44_10} we have, for all
$n \ge N$ and $t \ge 0$,
\be
\label{p44_11}
{\frac d{dt}} \| \vb_{n,2}\|^2_{H^1}
+ {\frac 32} \nu \| \nabla \vb_{n,2} \|^2
\le 
c \epsilon  \left (    \|  v \|_{H^1}^4 + \| v_t \|_{H^1}^2 
+ |y(\theta_t \omega) |^2 
+ |y(\theta_t \omega) |^4 \right ).
  \ee
  Note that \eqref{deltacond} and \eqref{p44_10} imply
  $$ {\frac 32} \nu \| \nabla \vb_{n,2} \|^2
  \ge
   \nu \| \nabla \vb_{n,2} \|^2
   + {\frac 12} \nu \lambda_{n+1} \| \vb_{n,2} \|^2
   \ge
   \nu \| \nabla \vb_{n,2} \|^2
   + {\frac 12} \nu \lambda  \| \vb_{n,2} \|^2
   \ge \delta \| \vb_{n,2} \|^2_{H^1},
   $$
   which along with \eqref{p44_11} shows that,
   for all $n\ge N$ and $t\ge 0$,
   \be
\label{p44_12}
{\frac d{dt}} \| \vb_{n,2}\|^2_{H^1}
+ \delta  \|   \vb_{n,2} \|^2_{H^1}
\le 
c \epsilon  \left (    \|  v \|_{H^1}^4 + \| v_t \|_{H^1}^2 
+ |y(\theta_t \omega) |^2 
+ |y(\theta_t \omega) |^4 \right ).
  \ee
  Integrating  \eqref{p44_12} over $(0,t)$, we find that, for all
  $n\ge N$ and $t\ge 0$,
  $$\| \vb_{n,2} (t, \omega  )\|^2_{H^1}
  \le e^{-\delta t } \| \vb_{n,2} (0, \omega  )\|^2_{H^1}
  $$
$$
  + c \epsilon
  \int_0^t e^{\delta (s-t)}
  \left (
  1 + \| v(s, \omega, v_0(\omega) )\|^4_{H^1}
  + \| v_s (s, \omega, v_0(\omega ) ) \|^2_{H^1}
  + |y(\theta_{s} \omega) |^2 
+ |y(\theta_{s}\omega) |^4
  \right ) ds.
 $$
Replacing $\omega$ by $\theta_{-t} \omega$ in the above, we get that,
for all
  $n\ge N$ and $t\ge 0$,
  $$\| \vb_{n,2} (t, \theta_{-t} \omega  )\|^2_{H^1}
  \le e^{-\delta t } \| \vb_{n,2} (0, \theta_{-t} \omega  )\|^2_{H^1}
  $$
   \be
  \label{p44_20}
  + c \epsilon
  \int_0^t e^{\delta (s-t)}
  \left (
  1 + \| v(s, \theta_{-t}\omega, v_0(\theta_{-t}\omega) )\|^4_{H^1}
  + \| v_s (s, \theta_{-t}\omega, v_0(\theta_{-t}\omega ) ) \|^2_{H^1}
  + |y(\theta_{s-t} \omega) |^2 
+ |y(\theta_{s-t}\omega) |^4
  \right ) ds.
\ee
For the first term on the right-hand side of
\eqref{p44_20}, by \eqref{vbv} we have
$$
e^{-\delta t } \| \vb_{n,2} (0, \theta_{-t} \omega  )\|^2_{H^1}
\le 
e^{-\delta t } \| \vb_0 (  \theta_{-t}\omega  )\|^2_{H^1}
 \le 
c e^{-\delta t } \| v_0 (  \theta_{-t}\omega  )\|^2_{H^1},
$$
which converges to zero as $t \to \infty$, and hence there is
$T_1 =T_1(B, \omega, \epsilon)>0$ such that for all
$t\ge T_1$,
\be
\label{p44_21}
e^{-\delta t } \| \vb_{n,2} (0, \theta_{-t} \omega  )\|^2_{H^1}
\le  \epsilon.
\ee
For the second term on the right-hand side 
of
\eqref{p44_20}, by \eqref{vs}, Lemma \ref{lem42} and
the proof of Lemma \ref{lem43} we can show that there are
$T_2=T_2 (B, \omega, \epsilon)>0$
and $r(\omega)>0$  such that for all
$t\ge T_2$,
 $$  
  \int_0^t e^{\delta (s-t)}
  \left (
  1 + \| v(s, \theta_{-t}\omega, v_0(\theta_{-t}\omega) )\|^4_{H^1}
  + \| v_s (s, \theta_{-t}\omega, v_0(\theta_{-t}\omega ) ) \|^2_{H^1}
  + |y(\theta_{s-t} \omega) |^2 
+ |y(\theta_{s-t}\omega) |^4
  \right ) ds
  $$
 \be
 \label{p44_30}
  \le c  (1+ r(\omega) ).
\ee
The details for the proof of  \eqref{p44_30}  are omitted here.
Let $T =\max\{T_1, T_2\}$.
Then by 
\eqref{p44_20}-\eqref{p44_30} we  get that, for all
$n\ge N$ and $t \ge T$,
$$\| \vb_{n,2} (t, \theta_{-t} \omega  )\|^2_{H^1}
\le \epsilon + c\epsilon (1 + r(\omega) ),
$$
which completes the proof.
 \end{proof}

\begin{lem}
\label{lem45}
 Assume that $g \in \h$, $ h  \in  H^1_0 (Q)$ and
 \eqref{fcond} holds. Let
 $B=\{B(\omega)\}_{\omega \in \Omega}\in \mathcal{D}$, 
 $t_m \to \infty$ and $v_{0,m}  \in B(\theta_{-t_m} \omega )$.
 Suppose $v (t, \theta_{-t} \omega, v_{0,m} )$ satisfies 
   \eqref{ve}-\eqref{ve2} with initial condition
 $v_{0,m}$ and
 $$
 \vb _m (x, t, \theta_{-t} \omega ) = \psi \left (
 {\frac {x_3^2}{k^2}}
 \right ) v   (x, t, \theta_{-t} \omega, v_{0,m} ),
 $$ where $k\ge 1$ is fixed. Then
 for $P$-a.e. $\omega \in \Omega$, 
 the sequence $\{\vb_m (t_m, \theta_{-t_m} \omega )\}_{m=1}^\infty$
 has a  convergent  subsequence  in $H^1_0(Q)$.
 \end{lem}
 \begin{proof}
 By Lemma \ref{lem41}, for $P$-a.e.  $\omega \in \Omega$, 
 there is $T_1 =T_1(B, \omega)$ such that for all $t \ge T_1$,
 \be
 \label{p45_1}
 \| v(t, \theta_{-t} \omega, v_0(\theta_{-t} \omega ) ) \|_{H^1_0(Q)}
 \le r(\omega),
 \ee
 where  $r(\omega)$ is a positive random function.
 Since $t_m \to \infty$, there is $M_1 =M_1(B, \omega)$ such that for all
 $m \ge M_1$,
 $t_m \ge T_1$, and hence by \eqref{p45_1} we have
$$
 \| v(t_m, \theta_{-t_m} \omega, v_0(\theta_{-t_m} \omega ) ) \|_{H^1_0(Q)}
 \le r(\omega), \quad \forall \ m \ge M_1,
$$
which implies that
\be
\label{p45_2}
 \| v(t_m, \theta_{-t_m} \omega, v_{0,m} ) \|_{H^1_0(Q)}
 \le r(\omega), \quad \forall \ m \ge M_1.
\ee
By \eqref{p45_2} we find that
\be
\label{p45_3}
 \| \vb_m  (t_m, \theta_{-t_m} \omega  ) \|_{H^1_0(Q)}
 \le c   r(\omega), \quad \forall \ m \ge M_1.
\ee
Given $\epsilon>0$, it follows from Lemma \ref{lem44} that
there are $T_2 =T_2(B, \omega, \epsilon)$ and
$N  =N (\omega, \epsilon)$ such that for all
$t\ge T_2$,
\be
\label{p45_4}
 \| (I-P_N) \vb (t , \theta_{-t } \omega  ) \|_{H^1_0(Q_{2k})}
 \le \epsilon.
\ee
Take $M_2 =M_2(B, \omega, \epsilon)$ large enough such that
$t_m \ge T_2$ for $m \ge M_2$. Then we get from \eqref{p45_4}
that
\be
\label{p45_5}
 \| (I-P_N) \vb _m (t_m , \theta_{-t_m } \omega  ) \|_{H^1_0(Q_{2k})}
 \le \epsilon, \quad \forall \ m \ge M_2.
\ee
On the other hand, \eqref{p45_3} shows that
the sequence $\{P_N \vb_m  (t_m, \theta_{-t_m} \omega  )\}$
is bounded in the finite-dimensional space 
$P_N H^1_0(Q_{2k})$ and hence is precompact
in  $P_N H^1_0(Q_{2k})$, which along with
\eqref{p45_5} 
implies  the precompactness of
 $\{\vb_m (t_m , \theta_{-t_m } \omega  )\}$
 in $H^1_0(Q_{2k})$.
 Note that $\vb_m (x, t_m , \theta_{-t_m } \omega  )
 =0$ for $x \notin Q_{2k}$
 and hence  
 $\{\vb_m (t_m , \theta_{-t_m } \omega  )\}$
 is precompact in $H^1_0(Q)$.
 \end{proof}
  
 Next we establish the asymptotic compactness of
 the solutions of problem \eqref{ve}-\eqref{ve3}.
\begin{lem}
\label{lem46}
 Assume that $g \in \h$, $ h  \in  H^1_0(Q)$ and
 \eqref{fcond} holds. Let
 $B=\{B(\omega)\}_{\omega \in \Omega}\in \mathcal{D}$, 
 $t_n \to \infty$ and $v_{0,n}  \in B(\theta_{-t_n} \omega )$.
  Then
 for $P$-a.e. $\omega \in \Omega$, 
 the sequence $\{v (t_n, \theta_{-t_n} \omega, v_{0,n} )\}_{n=1}^\infty$
 has a  convergent  subsequence  in $H^1_0(Q)$.
 \end{lem}
 
\begin{proof}
Given $\epsilon>0$, it follows from Lemma \ref{lem43} that, 
   for $P$-a.e.  $\omega \in \Omega$, 
 there are $T_1 =T_1(B, \omega, \epsilon)$
and $k_0 =k_0 (\omega, \epsilon)$  such that for all $t \ge T_1$,
 \be
 \label{p46_1}
 \int_{Q\backslash Q_{k_0}} \left ( |v(t, \theta_{-t} \omega, v_0(\theta_{-t} \omega )
  )|^2
  +
   |\nabla v(t, \theta_{-t} \omega, v_0(\theta_{-t} \omega )
  )|^2 \right ) dx
  \le \epsilon.
 \ee
 Let $N_1=N_1(B, \omega, \epsilon)$ be large enough such that
 $t_n \ge T_1$ for $n \ge N_1$. Then by \eqref{p46_1} we have
$$
 \int_{Q\backslash Q_{k_0}} 
\left ( |v(t_n, \theta_{-t_n} \omega, v_0(\theta_{-t_n} \omega )
  )|^2
  +
   |\nabla v(t_n, \theta_{-t_n} \omega, v_0(\theta_{-t_n} \omega )
  )|^2 \right ) dx
  \le \epsilon, \quad \forall \ n \ge N_1,
$$
 which implies that
  \be
 \label{p46_2}
 \int_{Q\backslash Q_{k_0}} 
\left ( |v(t_n, \theta_{-t_n} \omega, v_{0,n}
  )|^2
  +
   |\nabla v(t_n, \theta_{-t_n} \omega, v_{0,n}
  )|^2 \right ) dx
  \le \epsilon, \quad \forall \ n \ge N_1.
  \ee
  Denote by
  $$
 \vb _n (x, t_n, \theta_{-t_n} \omega ) = \psi \left (
 {\frac {x_3^2}{k^2}}
 \right ) v   (x, t_n, \theta_{-t_n} \omega, v_{0,n} ).
 $$
 Then from Lemma \ref{lem45} we know that,
 up to a subsequence, $\{\vb _n ( t_n, \theta_{-t_n} \omega ) \}$
 is convergent  in $H^1_0(Q)$, which shows that
 $\{\vb _n ( t_n, \theta_{-t_n} \omega ) \}$
 is a Cauchy sequence  in $H^1_0(Q)$, and hence
 also  a  Cauchy sequence  in $H^1(Q_{k_0})$.
 Note that  $ \vb _n ( t_n, \theta_{-t_n} \omega )  $
 $ =  v  ( t_n, \theta_{-t_n} \omega , v_{0,n} )  $
 in $Q_{k_0}$ and thus
 $ \{ v  ( t_n, \theta_{-t_n} \omega , v_{0,n} ) \} $
 is a Cauchy sequence in $H^1(Q_{k_0})$.
 This along with \eqref{p46_2} shows that
 $ \{ v  ( t_n, \theta_{-t_n} \omega , v_{0,n} ) \} $
 is a Cauchy sequence in $H^1_0(Q)$, as desired.
\end{proof}

\section{Random attractors
  }
\setcounter{equation}{0}

In this section, we prove the existence of a 
$\mathcal{D}$-random attractor  for the 
random dynamical system $\Phi$ associated with
the stochastic \bbme on the unbounded channel $Q$.
By \eqref{uv} and \eqref{phi}, $\Phi$ satisfies
\be
\label{att1}
\Phi ( t, \theta_{-t} \omega, u_0(\theta_{-t} \omega ) )
=u(t, \theta_{-t} \omega, u_0(\theta_{-t} \omega ) )
= v(t, \theta_{-t} \omega, v_0(\theta_{-t} \omega ) ) + z(\omega),
\ee
where $v_0(\theta_{-t} \omega )
= u_0(\theta_{-t} \omega )
- z(\theta_{-t} \omega)$.
Let $B=\{B(\omega)\}_{\omega \in \Omega} \in \mathcal{D}$ and define
\be
\label{att2}
{\tilde{B}}(\omega)
= \{ v \in H^1_0(Q): \ 
\| v \|_{H^1_0} \le \| u(\omega) \|_{H^1_0} + \| z(\omega) \|_{H^1_0} ,
\quad u(\omega)  \in B(\omega)
\}.
\ee
We  claim that ${\tilde{B}} =\{ {\tilde{B}}(\omega)\}_{\omega \in
\Omega} $ belongs to $  \mathcal{D}$ provided 
$B=\{B(\omega)\}_{\omega \in \Omega} \in \mathcal{D}$.
Note that $B=\{B(\omega)\}_{\omega \in \Omega} \in \mathcal{D}$
implies that
\be
\label{att3}
\lim_{t \to \infty} e^{-{\frac 18} \delta t }  d (B( \theta_{-t} \omega  ) )
=0.
\ee
Since $z(\omega)$ is tempered, by \eqref{att2}-\eqref{att3}  we have
$$
\lim_{t \to \infty} e^{-{\frac 18} \delta t }d({\tilde{B}} (\theta_{-t} \omega ) )
\le \lim_{t \to \infty} e^{-{\frac 18} \delta t }d({{B}} (\theta_{-t} \omega ) )
+ \lim_{t \to \infty} e^{-{\frac 18} \delta t }\| z (\theta_{-t} \omega )\|_{H^1_0}
=0,
$$
which shows 
${\tilde{B}} =\{ {\tilde{B}}(\omega)\}_{\omega \in
\Omega}  \in  \mathcal{D}$.
Then by Lemma \ref{lem41}, for $P$-a.e. $\omega \in \Omega$,
if $v_0(\omega) \in {\tilde{B}}(\omega)$, there is $T_1 =T_1({\tilde{B}}, \omega)
$ such that   for all $t \ge T_1$,
\be
\label{att8}
\| v(t, \theta_{-t} \omega, v_0(\theta_{-t} \omega ) ) \|_{H^1_0}
\le r(\omega ),
\ee
where $r(\omega)$ is a positive random function satisfying
\be
\label{att9}
e^{-{\frac 18} \delta t} r(\theta_{-t} \omega )
\to 0 \quad 
\mbox{as} \ t \to \infty.
\ee
Denote by
\be
\label{att10}
K(\omega) = \{ u \in H^1_0(Q): \
\| u \|_{H^1_0} 
\le r(\omega)  + \| z(\omega ) \|_{H^1_0} \}.
\ee
Then by \eqref{att9}  we have
$$
\lim_{t \to \infty} e^{-{\frac 18} \delta t }d(K(\theta_{-t} \omega ) )
\le \lim_{t \to \infty} e^{-{\frac 18} \delta t } r (\theta_{-t} \omega )  
+ \lim_{t \to \infty} e^{-{\frac 18} \delta t }\| z (\theta_{-t} \omega )\|_{H^1_0}
=0,
$$
which implies that $K =\{K(\omega)\}_{\omega \in \Omega} 
\in \mathcal{D}$. We now show that
$K$ is also an absorbing  set of $\Phi$ in
$\mathcal{D}$.
Given  $B=\{B(\omega)\}_{\omega \in \Omega } \in \mathcal{D}$ 
and $u_0(\omega) \in B(\omega)$, 
by \eqref{att1} and \eqref{att8} we find that, for all
$t \ge T_1$,
$$
\| u(t, \theta_{-t} \omega, u_0(\theta_{-t} \omega ) ) \|_{H^1_0}
\le 
\| v(t, \theta_{-t} \omega, v_0(\theta_{-t} \omega ) ) \|_{H^1_0}
+ 
\| z  (\omega)   \|_{H^1_0}
\le r(\omega) +  \| z ( \omega )  \|_{H^1_0},
$$
which along with \eqref{att1} and  \eqref{att10} implies  that
\be
\label{att11}
\Phi (   t, \theta_{-t} \omega, B(\theta_{-t} \omega ) )
\subseteq K(\omega),\quad \forall \ t \ge T_1,
\ee
and hence $K =\{K(\omega)\}_{\omega \in \Omega } \in \mathcal{D}$ is 
a closed  absorbing  set of $\Phi$ in
$\mathcal{D}$. In a word, we have proved the following result.

\begin{lem}
\label{lem51}
 Assume that $g \in \h$, $ h  \in  H^1_0 (Q)$ and
 \eqref{fcond} holds. Let
 $K=\{K(\omega)\}_{\omega \in \Omega}$ be given by
\eqref{att10}. Then   
$K =\{K(\omega)\}_{\omega \in \Omega } \in \mathcal{D}$ is 
a closed  absorbing  set of $\Phi$ in
$\mathcal{D}$.
 \end{lem}
 
 As an immediate consequence of Lemma \ref{lem46}, we find that
 $\Phi$ is $\mathcal{D}$-pullback  asymptotically compact
 in $H^1_0(Q)$.

\begin{lem}
\label{lem52}
 Assume that $g \in \h$, $ h  \in  H^1_0 (Q)$ and
 \eqref{fcond} holds. 
 Then the random dynamical system $\Phi$ is
 $\mathcal{D}$-pullback  asymptotically compact
 in $H^1_0(Q)$; that is,
 for $P$-a.e. $\omega \in \Omega$, the sequence 
 $\{\Phi (t_n, \theta_{-t_n} \omega, u_{0,n}  )\}$ has a 
 convergent  subsequence in $\hone$ provided 
 $t_n \to \infty$, 
$B =\{B(\omega)\}_{\omega \in \Omega } \in \mathcal{D}$
and $u_{0,n} \in B(\theta_{-t_n} \omega )$.
 \end{lem}
 
 \begin{proof}
 Since  $B =\{B(\omega)\}_{\omega \in \Omega }$
belongs to  $ \mathcal{D}$, so does 
${\tilde{B}} =\{{\tilde{B}}(\omega)\}_{\omega \in \Omega }$
which is given by \eqref{att2}.
Then it follows from Lemma \ref{lem46} that, 
for $P$-a.e. $\omega \in \Omega$, up to a subsequence,
$\{v(t_n, \theta_{-t_n} \omega, v_{0,n} )\}$ is convergent 
in $\hone$, where
$v_{0,n} = u_{0,n} - z(\theta_{-t_n } \omega )$
$\in {\tilde{B}}(\theta_{-t_n} \omega )$.
This along with  \eqref{att1} shows that, up to a subsequence,
$\{\Phi (t_n, \theta_{-t_n} \omega, u_{0,n}  )\}$  is  
 convergent in $\hone$.
 \end{proof}
 
 We are now in   a position to establish the existence 
 of a $\mathcal{D}$-random attractor for $\Phi$.

\begin{thm}
\label{thm51}
 Assume that $g \in \h$, $ h  \in  H^1_0 (Q)$ and
 \eqref{fcond} holds.  Let $\mathcal{D}$ be the collection
of random sets given by
\eqref{Dset}.  Then the random dynamical 
 system $\Phi$ has a unique $\mathcal{D}$-random
 attractor in $H^1_0(Q)$.
 \end{thm}

\begin{proof}
Notice that $\Phi$ has a closed  absorbing set
$\{K(\omega)\}_{\omega \in \Omega}$ in $\mathcal{D}$ by Lemma
\ref{lem51}, and is $\mathcal{D}$-pullback asymptotically compact
in $\hone$ by Lemma \ref{lem52}. Hence the existence of a unique
$\mathcal{D}$-random attractor  follows from Proposition
\ref{att}
 immediately.
\end{proof}


\begin{thebibliography}{99}



\bibitem{ant1}
F. Antoci and M. Prizzi, Reaction-Diffusion equations on unbounded
thin domains,
{\em Topological Methods in Nonlinear  Analysis},
{\bf 18} (2001), 283-302.


\bibitem{ant2}
F. Antoci and M. Prizzi, Attractors and global averaging of non-autonomous
Reaction-Diffusion equations in $\mathbb{R}^n$.
{\em Topological Methods in Nonlinear  Analysis},
{\bf 20} (2002), 229-259.


\bibitem{arn1}
L. Arnold, {\em Random Dynamical Systems}, Springer-Verlag, 1998.

\bibitem{arr1}
J. M. Arrieta, J. W. Cholewa, T. Dlotko,  and A. Rodriguez-Bernal,
Asymptotic behavior and attractors for Reaction Diffusion equations
in unbounded  domains,
{\em Nonlinear Analysis}, {\bf 56} (2004), 515-554.



  \bibitem{ast-bis-bis-fer}
 M.A. Astaburuaga, E. Bisognin,
V. Bisognin,  and  C. Fernandez,
 Global attractors and finite dimensionality for
 a class of dissipative equations of BBM's type,
  {\em Electron. J. Differential Equations},   {\bf 25}  (1998),
   1-14.


  \bibitem{avr}
  J. Avrin, The generalized \bbme in $\rn$ with
  singular initial data,
  {\em  Nonlinear Anal.},  {\bf 11}  (1987), 139-147.


 \bibitem{avr-gol}
 J. Avrin and  J.A. Goldstein, Global existence for the
 \bbme in arbitrary dimensions,
 {\em Nonlinear Anal.},  {\bf 9}  (1985), 861-865.



\bibitem{bab1}
A.V. Babin and  M.I.  Vishik,     Attractors of
Evolution Equations,  North-Holland, Amsterdam,
1992.


\bibitem{bal1}
J.M.   Ball,  Continuity properties and global
attractors of generalized semiflows and the Navier-Stokes
equations, {\em  J. Nonl. Sci.},  {\bf 7} (1997),
475-502.


\bibitem{bal2}
J.M.   Ball,
Global attractors for damped semilinear wave equations,
  {\em  Discrete Contin.  Dyn.  Syst.}   {\bf 10}   (2004),
31-52.



\bibitem{bat1}
P.W. Bates,  H. Lisei and  K.  Lu,
 Attractors for stochastic lattice
 dynamical systems,
{\em Stoch. Dyn.},   {\bf 6}  (2006),      1-21.



\bibitem{bat2}
P.W. Bates,   K.  Lu   and B. Wang, 
 Random Attractors for  Stochastic Reaction-Diffusion Equations
on Unbounded Domains,  submitted.



 \bibitem{ben-bon-mah} T.B. Benjamin,  J.L. Bona and  J.J. Mahony,
 Model equations  for long waves
 in  nonlinear  dispersive systems,
{\em Philos. Trans. Roy. Soc. London},  {\bf 272} (1972),  47-78.


  \bibitem{bon-bry} J.L. Bona and  P.J. Bryant, A mathematical
 model for long waves generated by wavemakers in nonlinear
 dispersive
 systems,
 {\em Proc. Cambridge Philos.
 Soc.},  {\bf 73} (1973), 391-405.


\bibitem{bon-dou} J.L. Bona and  V.A. Dougalis, An
 initial - and boundary-value
 problem  for a model equation
 for propagation of long waves,
 {\em J. Math. Anal. Appl.},
  {\bf 75} (1980), 503-522.



    \bibitem{car1}
    T. Caraballo, J. A. Langa and J. C. Robinson,
    a stochastic pitchfork bifurcation in a reaction-diffusion equation,
    {\em  Proc. R. Soc. Lond. A}, {\bf 457} (2001),  2041-2061.




 \bibitem{cel-kal-pol}
  A.O. Celebi, V.K. Kalantarov and  M. polat,
  Attractors for the generalized \bbm equation,
  {\em J. Differential Equations},  {\bf 157}  (1999), 439-451.

  \bibitem{che}
  Y.M. Chen, Remark on the global existence
  for the generalized \bbm equations in arbitrary dimension,
  {\em Appl. Anal.},  {\bf 30}  (1988), 1-15.


  \bibitem{chu-pol-sie}
  I. Chueshov, M. Polat and  S. Siegmund,  Gevrey regularity of global attractor
  for generalized \bbm equation,
   {\em Mat. Fiz. Anal. Geom.},   {\bf 11}   (2004),   226-242.




    \bibitem{cra1}
    H. Crauel,  A. Debussche and
    F. Flandoli, Random attractors,
    {\em J. Dyn. Diff. Eqns.}, {\bf 9} (1997), 307-341.


    \bibitem{cra2}
    H. Crauel  and
    F. Flandoli,  Attractors for random dynamical systems,
    {\em Probab. Th. Re. Fields}, {\bf 100} (1994), 365-393.

     \bibitem{fla1}
    F. Flandoli and B. Schmalfu$\beta$,
Random attractors for
the 3D stochastic Navier-Stokes equation with multiplicative
noise,
    {\em Stoch. Stoch. Rep.}, {\bf 59} (1996),  21-45.



\bibitem{ghi1}
J.M. Ghidaglia, A note on the strong convergence towards attractors for damped
forced KdV equations,
{\em J. Differential Equations}, {\bf 110} (1994), 356-359.


 
 
   \bibitem{gol-wic}
   J.A. Goldstein and B.J.  Wichnoski, On the \bbme
   in higher dimensions,
{\em Nonlinear Anal.},  {\bf 4}  (1980),  665-675.



\bibitem{gou1} O. Goubet and  R. Rosa,
Asymptotic smoothing and the global attractor of a weakly
damped KdV equation on the real line,
 {\em J. Differential Equations},   {\bf 185}  (2002),  25-53.


 \bibitem{hal1}
  J.K.  Hale,   Asymptotic Behavior of
Dissipative Systems,
American Mathematical Society,
  Providence, RI, 1988.

\bibitem{ju1}
N. Ju,  The $H^1$-compact global attractor for the solutions
to the Navier-Stokes equations in two-dimensional unbounded domains,
  {\em  Nonlinearity},  {\bf 13} (2000),
1227-1238.



 \bibitem{med-mil} L.A. Medeiros,  M. Milla Miranda,
 Weak  solutions  for a nonlinear dispersive  equation,
  {\em J. Math. Anal. Appl.}, {\bf 59} (1977), 432-441.

 \bibitem{med-per} L.A. Medeiros, G. Perla Menzala,
 Existence  and  uniqueness   for periodic
 solutions  of the
 Benjamin-Bona-Mahony equation,
 {\em SIAM J. Math. Anal.},  {\bf 8} (1977), 792-799.
 
 \bibitem{mor1}
 F. Morillas and J. Valero,
 Attractors for Reaction-Diffusion equations in $\mathbb{R}^n$ with
 continuous  nonlinearity, {\em Asymptotic Analysis}, {\bf 44} (2005),
 111-130.

\bibitem{moi2}
I. Moise and R. Rosa, on the regularity of the global attractor of a weakly damped,
forced Korteweg-de Vries equation,
{\em Adv. Differential  Equations}, {\bf 2} (1997), 257-296.

\bibitem{moi1}
I. Moise, R. Rosa and X. Wang,
Attractors for non-compact semigroups via energy equations,
  {\em  Nonlinearity},  {\bf 11} (1998),  1369-1393.

\bibitem{pri1}
M. Prizzi, Averaging, Conley index continuation and recurrent
dynamics in almost-periodic equations,
{\em J. Differential Equations}, {\bf 210} (2005), 429-451.

    \bibitem{rob1}
 J.C. Robinson, Infinite-Dimensional Dynamical Systems,
 Cambridge University Press, Cambridge, UK, 2001.

\bibitem{rod1}
A. Rodrigue-Bernal and B. Wang,
Attractors for partly dissipative Reaction Diffusion systems
in $\mathbb{R}^n$,
{\em J. Math. Anal. Appl.}, {\bf 252} (2000),  790-803.

 \bibitem{ros1}
R. Rosa,  The global attractor for the 2D Navier-Stokes flow on some
unbounded  domains,
  {\em  Nonlinear Anal.},  {\bf 32} (1998),   71-85.

    \bibitem{sel1}
 R. Sell and  Y. You,
 Dynamics of Evolutionary Equations,
 Springer-Verlag,
New York, 2002.

\bibitem{sta1}
M. Stanislavova, A. Stefanov and B. Wang,
Asymptotic smoothing and attractors for the generalized
 Benjamin-Bona-Mahony equation
on $\R^3$,
{\em J. Differential Equations}, {\bf 219} (2005), 451-483.

 \bibitem{sun1}
 C. Sun and C. Zhong, Attractors for the semilinear Reaction-Diffusion
 equation with distribution derivatives in unbounded  domains,
 {\em Nonlinear Analysis}, {\bf 63} (2005), 49-65.

\bibitem{tem1}
 R.  Temam,    Infinite-Dimensional Dynamical
Systems in Mechanics and Physics,
  Springer-Verlag,
  New York, 1997.



\bibitem{wan1}
B. Wang,  Attractors for reaction-diffusion equations in unbounded domains,
{\em Physica D}, {\bf 128}  (1999), 41-52.


 
 \bibitem{wan2}
 B. Wang,
 Regularity of attractors for the
 \bbm equation,
  {\em J. Phys. A, Math. Gen.},  {\bf 31}  (1998), 7635-7645.

 \bibitem{wan3}
 B. Wang, Strong Attracrors for
 the Benjamin-Bona-Mahony equation,
 {\em Appl. Math. Lett.},  {\bf 10}  (1997),  23-28.


 \bibitem{wan-yan}
 B. Wang, W. Yang, Finite dimensional behaviour
 for the  Benjamin-Bona-Mahony equation,
 {\em J. Phys. A, Math. Gen.},  {\bf 30} (1997), 4877-4885.


\bibitem{wanx}
X. Wang, An energy equation for the weakly damped driven nonlinear Schrodinger equations
and its applications, {\em Physica D}, {\bf  88} (1995),  167-175.

 
  

 \end{thebibliography}
 \end{document}